\documentclass{scrartcl}
\usepackage{amsmath}
\usepackage{amsfonts}
\usepackage{amssymb}
\usepackage{dsfont}
\usepackage[T1]{fontenc}
\usepackage[latin1]{inputenc}
\usepackage{epsfig}
\usepackage{graphicx}

\usepackage[boxed, lined]{algorithm2e}
\usepackage{float}
\usepackage{comment}
\usepackage{mathtools}

\newcommand{\ba}{\mathbf{a}}
\newcommand{\bb}{\mathbf{b}}

\newcommand{\by}{\mathbf{y}}
\newcommand{\bw}{\mathbf{w}}
\newcommand{\bv}{\mathbf{v}}

\newcommand{\D}{{\mathcal{D}}}

\newcommand{\bA}{{\mathbf{A}}}
\newcommand{\bB}{{\mathbf{B}}}
\newcommand{\Nu}{{\mathcal{N}}}

\newcommand{\N}{\mathbb{N}}

\newcommand{\R}{\mathbb{R}}

\newcommand{\Rd}{\mathbb{R}^d}

\newcommand{\bE}{\mathbf{E}}

\newcommand{\beq}{\begin{eqnarray*}}

\newcommand{\eeq}{\end{eqnarray*}}

\newcommand{\beqm}{\begin{eqnarray}}

\newcommand{\eeqm}{\end{eqnarray}}

\newtheorem{theorem}{Theorem}

\newtheorem{lemma}{Lemma}

\newcommand{\EXP}{{\mathbf E}}
\newcommand{\PROB}{{\mathbf P}}

\renewcommand{\bf}{\normalfont \bfseries}
\renewcommand{\it}{\normalfont \itshape}

\allowdisplaybreaks
\begin{document}
\renewcommand{\thefootnote}{\fnsymbol{footnote}}
\newcommand{\F}{{\cal F}}
\newcommand{\Sp}{{\cal S}}
\newcommand{\G}{{\cal G}}
\newcommand{\HH}{{\cal H}}
%\maketitle
%\noindent

\begin{center}

  {\LARGE \bf
    On the universal consistency of an over-parametrized
    deep neural network estimate learned by gradient descent
  }
\footnote{
Running title: {\it Consistent neural networks}}
\vspace{0.5cm}

Selina Drews\footnote{Corresponding author. Tel:
  +49-6151-16-23372, Fax:+49-6151-16-23381}
 and
Michael Kohler

{\it 
Fachbereich Mathematik, Technische Universit\"at Darmstadt,
Schlossgartenstr. 7, 64289 Darmstadt, Germany,
email: drews@mathematik.tu-darmstadt.de, kohler@mathematik.tu-darmstadt.de}

\end{center}
\vspace{0.5cm}

\begin{center}
August 22, 2022
\end{center}
\vspace{0.5cm}

\noindent
    {\bf Abstract}\\
    Estimation of a multivariate regression function from independent
    and identically distributed data is considered. An estimate
    is defined which fits a deep neural network consisting of
    a large number of fully connected neural networks, which
    are computed in parallel, via gradient descent to the data. The estimate
    is over-parametrized in the sense that the number
    of its parameters is much larger than the sample size. It is shown
    that in case of a suitable random initialization of the network,
    a suitable
    small stepsize of the gradient descent, and a number of gradient
    descent steps which is slightly larger than the reciprocal of the
    stepsize of the gradient descent, the estimate is universally
    consistent in the sense that its expected $L_2$ error converges
    to zero for all distributions of the data where the response
    variable is square integrable.
    
    \vspace*{0.2cm}


\vspace*{0.2cm}


\section{Introduction}
\label{se1}
Deep neural networks belong nowadays to the most promising approaches
in many different applications. They have been successfully
applied, e.g., in
image classification
(cf., e.g., Krizhevsky, Sutskever and Hinton  (2012)),
text classification
(cf., e.g., Kim (2014)),
machine translation
(cf., e.g., Wu et al. (2016)) 
or
mastering of games
(cf., e.g., Silver et al. (2017)).

In the last few years
various results concerning
the approximation power of deep neural networks (cf., e.g., 
Yarotsky (2017),  Yarotsky and Zhevnerchuck (2020),
Lu et al.  (2020), Langer (2021b)
and the literature cited 
therein)
or concerning the statistical risk of
corresponding least squares estimates
(cf., e.g.,
Bauer and Kohler (2019), Kohler and Krzy\.zak (2017),
Schmidt-Hieber (2020), 
Kohler and Langer (2021), 
Langer (2021a), Imaizumi and Fukumizu (2019), Suzuki (2018), Suzuki 
and Nitanda (2019),  and the literature cited 
therein)
have been derived.

The above results ignore two important features of the typical
application of deep neural networks: Firstly, in practice the
estimates are computed using gradient descent and not (as in the
theoretical results above) the principle of least squares. And secondly,
often the applied neural networks are over-parametrized in the
sense that the number of parameters is much larger than their sample
size.

These two principles contradict classical theory. Nevertheless, to some surprise, they work very well in practice.
In Bartlett, Montanari and  Rakhlin (2021), the theory of this observation is explored in more detail. To do this, they put forward two hypotheses.
The first hypothesis concerns the \textit{tractability via over-parametrization}. It is conjectured that even if the objective function is non-convex, the hardness of the optimization problem depends on the relationship between the dimension of the parameter space (the number of optimization variables) and the sample size. Thus the tractability is given if and only if a model is chosen that is over-parametrized. This is in contrast to the classical assumption that statistical learning is achieved by restricting to linearly parameterized classes of functions and convex objectives. 
The second hypothesis concerns \textit{generalization via implicit regularization}.
In classical theory, one wanted to avoid over-parametrized neural networks and therefore restricted them to an under-parametrized regime or a suitable regularizing regime. It was assumed that a method that has too many degrees of freedom by perfectly interpolating noisy data cannot have a good generalization.
However, it was observed in practice that over-parametrized models generalize well. This is very interesting since empirical evidence shows that an optimization task is simplified if the model is sufficiently over-parametrized.

Bartlett, Montanari and  Rakhlin (2021) suggest that deep learning models can be divided into a simple component and a spiky component. The simple component is useful for prediction, and the spiky component is useful for overfitting.
If the model is suitably over-parametrized, this interpolation does not affect the prediction accuracy.

Nonconvex empirical risk minimization problems in a linear regime are solved efficiently by gradient methods.
In a linear regime, a parameterized function can be approximated exactly by its linearization over an initial parameter vector. Bartlett, Montanari
and  Rakhlin (2021) were able to show that for a suitable parameterization and initialization, a gradient method remains in the linear regime. Further, it leads to linear convergence of the empirical risk and to a solution whose prediction is well approximated by the linearization of the initialization. Especially, they showed that for two-layer networks in the linear regime, a suitably large over-parametrization together with a suitable initialization is sufficient.
This theory is not able to capture training schemes in which the weights genuinely change.

One approach beyond the linear regime considered in Bartlett, Montanari
and  Rakhlin (2021) is the mean field limit.  Here, the weights move in a nontrivial way during training, even though the network is infinitely wide.
Using mean field theory, global convergence results can be proved for two-layer (Mei, Montanari, and Nguyen (2018), Chizat and Bach (2018)) and multi-layer neural networks (Nguyen and Pham (2020)).

Another approach is to consider the linearized evolution as a Taylor expansion of the first order, and then construct higher-order approximations.
Other approaches are considered in Dyer and Gur-Ar (2019) and Hanin and Nica (2019).

It is well-known that gradient descent leads to
a small empirical $L_2$ risk in over-parametrized neural networks,
see, e.g.,
Allen-Zhu, Li and Song (2019), Kawaguchi and Huang (2019)
and the literature cited therein.
However, such results are in general not useful for the proof
of the consistency of corresponding estimates, because it
was shown in Kohler and Krzy\.zak (2021) that any estimate which
interpolates the training data does not generalize well in a sense
that its error does not converge to zero for sample size tending to infinity
in case of a general design measure.

In the current paper, we analyze deep neural networks in the context
of nonparametric regression.
Here we consider
an $\Rd \times \R$--valued random
vector $(X,Y)$ with $\EXP Y^2 < \infty$,
where $X$ is the so--called observation vector
and $Y$ is the  so-called response variable. We assume that
a sample of $(X,Y)$, i.e.,  a data set
\begin{equation}
  \label{inteq1}
\D_n = \left\{
(X_1,Y_1), \ldots, (X_n,Y_n) 
\right\},
\end{equation}
where
$(X,Y)$, $(X_1,Y_1)$, \ldots, $(X_n,Y_n)$ are i.i.d., is available. We are searching for an estimator
\[
m_n(\cdot)=m_n(\cdot, \D_n):\Rd \rightarrow \R
\]
of the so--called regression function $m:\Rd \rightarrow \R$,
$m(x)=\EXP\{Y|X=x\}$ such that the so--called $L_2$ error
\[
\int |m_n(x)-m(x)|^2 {\PROB}_X (dx)
\]
is ``small'' (cf., e.g., Gy\"orfi et al. (2002)
for a systematic introduction to nonparametric regression and
a motivation for the $L_2$ error). 

In Section \ref{se2} we introduce an over-parametrized
deep neural network estimate, where the weights are learned
by gradient descent. Our main result is that in case we
initialize our starting weights randomly in a proper way, and
proceed with a suitable number of gradient descent steps
with a sufficiently small constant stepsize, this estimate $m_n$
is universally consistent in the sense that
\[
  \EXP \int | m_n(x)-m(x)|^2 \PROB_X (dx)
  \rightarrow 0
  \quad
  (n \rightarrow \infty)
  \]
  holds for every distribution of $(X,Y)$ with $\EXP Y^2 < \infty$.

  For many years it is well-known that universally consistent
  regression estimates exist, see Stone (1977) for the first result
  in this respect and Gy\"orfi et al. (2002) for an extensive overview
  of such results. So it is not surprising that deep neural networks
  have this property, too. However, our main result presents interesting
  aspects of the application of gradient descent to deep neural networks
  which are useful for
  proving such universal consistency: Firstly, the over-parametrization
  is useful in our result since it ensures that a finite subset
  of the initially chosen inner weights have nice properties. Secondly,
  due to the fact that we use a small stepsize, gradient descent
  applied to a properly regularized empirical $L_2$ risk
  is able to adjust the outer weights in an optimal way.
  And thirdly, due to the fact that the number of gradient descent
  steps is only slightly larger than the reciprocal of the stepsize, the
  inner weights do not change drastically during our learning.
  Altogether this enables our estimate to perform a kind
  of representation guessing instead of representation learning.

  In our proofs we use techniques that have been
  introduced in Braun et al. (2021)
  in the context of the analysis of gradient descent of neural networks
  with one hidden layer.
  These techniques have been also applied in Kohler and Krzy\.zak (2022)
  to analyze the performance of over-parametrized neural
  networks with one hidden layer. But in contrast to
  Kohler and Krzy\.zak (2022) we do not control the complexity
  of our estimate by using a strong regularization term. Instead
  we combine the techniques introduced in Braun et al. (2021)
  with the approach of Li, Gu and Ding (2021), which suggested to
  analyze the complexity of over-parametrized
  neural networks with metric entropy bounds.

  Throughout this paper we will use the following notation:
  The sets of natural numbers, real numbers and nonnegative real numbers
are denoted by $\N$, $\R$ and $\R_+$, respectively. For $z \in \R$, we denote
the smallest integer greater than or equal to $z$ by
$\lceil z \rceil$.
The Euclidean norm of $x \in \Rd$
is denoted by $\|x\|$.
For $f:\R^d \rightarrow \R$
\[
\|f\|_\infty = \sup_{x \in \R^d} |f(x)|
\]
is its supremum norm.

Let $\F$ be a set of functions $f:\Rd \rightarrow \R$,
let $x_1, \dots, x_n \in \Rd$, set $x_1^n=(x_1,\dots,x_n)$ and let
$p \geq 1$.
A finite collection $f_1, \dots, f_N:\Rd \rightarrow \R$
  is called an $L_p$ $\varepsilon$--cover of $\F$ on $x_1^n$
  if for any $f \in \F$ there exists  $i \in \{1, \dots, N\}$
  such that
  \[
  \left(
  \frac{1}{n} \sum_{k=1}^n |f(x_k)-f_i(x_k)|^p
  \right)^{1/p}< \varepsilon.
  \]
  The $L_p$ $\varepsilon$--covering number of $\F$ on $x_1^n$
  is the  size $N$ of the smallest $L_p$ $\varepsilon$--cover
  of $\F$ on $x_1^n$ and is denoted by $\Nu_p(\varepsilon,\F,x_1^n)$.

  If $A$ is a subset of $\Rd$ and $x \in \Rd$, then we set
  $1_A(x)=1$ if $x \in A$ and $1_A(x)=0$ otherwise.
For $z \in \R$ and $\beta>0$ we define
$T_\beta z = \max\{-\beta, \min\{\beta,z\}\}$. If $f:\R^d \rightarrow
\R$
is a function and $\F$ is a set of such functions, then we set
$
(T_{\beta} f)(x)=
T_{\beta} \left( f(x) \right)$ and
\[
T_{\beta} \F = \{ T_\beta f \, : \, f \in \F \}.
\]

  In Section \ref{se2} we define our estimate.
In Section \ref{se3} we present our main result concerning
the universal consistency of our deep neural network estimate
learned by gradient descent. The proof
of the main result is
given in Section \ref{se4}.

\section{Definiton of the estimate}
\label{se2}

Let
$\sigma(x)=1/(1+e^{-x})$
be the logistic squasher.
In the sequel we
use a network topology where
we compute a linear combination of $K_n$ fully connected
neural networks with $L$ layers and $r$ neurons per layer, i.e.,
we define our neural network as follows: We set
\begin{equation}\label{se2eq1}
f_\bw(x) = \sum_{j=1}^{K_n} w_{1,1,j}^{(L)} \cdot f_{j,1}^{(L)}(x) 
\end{equation}
for some $w_{1,1,1}^{(L)}, \dots, w_{1,1,K_n}^{(L)} \in \mathbb{R}$, where $f_{j,1}^{(L)}$ are recursively defined by
\begin{equation}
  \label{se2eq2}
f_{k,i}^{(l)}(x) = \sigma\left(\sum_{j=1}^{r} w_{k,i,j}^{(l-1)} \cdot f_{k,j}^{(l-1)}(x) + w_{k,i,0}^{(l-1)} \right)
\end{equation}
for some $w_{k,i,0}^{(l-1)}, \dots, w_{k,i, r}^{(l-1)} \in \mathbb{R}$
$(l=2, \dots, L)$
and
\begin{equation}
  \label{se2eq3}
f_{k,i}^{(1)}(x) = \sigma \left(\sum_{j=1}^d w_{k,i,j}^{(0)} \cdot x^{(j)} + w_{k,i,0}^{(0)} \right)
\end{equation}
for some $w_{k,i,0}^{(0)}, \dots, w_{k,i,d}^{(0)} \in \mathbb{R}$.

The above neural network consists of $K_n$ fully connected
neural networks with depth $L$, which are computed in parallel.
These networks have $r$ neurons in all layers except for the last layer,
where they only have one neuron. In the $k$-th such network we
denote the output of neuron $i$ in the $l$-th layer by $f_{k,i}^{(l)}$,
and the weight between neuron $j$ in the $(l-1)$-th layer and neuron
$i$ in the $l$-th layer is denoted by $w_{k,i,j}^{(l-1)}$. The number
of weights of the above neural network is given by
\[
K_n \cdot (1 + (L-1) \cdot r \cdot (r+1) + r \cdot (d+1)).
\]

In order to learn the
weight vector
$\bw=(w_{k,i,j}^{(l)})_{k,i,j,l}$
of our neural network we apply gradient descent to
a properly regularized empirical $L_2$ risk of our estimate.
We initialize
$\bw^{(0)}$ by setting
\begin{equation}
  \label{se2eq5}
  w_{1,1,j}^{(L)}=0 \quad \mbox{for } j=1, \dots, K_n,
\end{equation}
and by choosing all others weights
randomly such that all weights $w_{k,i,j}^{(l)}$ are independently
randomly chosen, and such that all weights
 $w_{k,i,j}^{(l)}$
with $1 \leq l<L$
are 
uniformly
distributed on $[-20 d \cdot (\log n)^2, 20 d \cdot (\log n)^2]$,
and all weights $w_{k,i,j}^{(0)}$
are 
uniformly
distributed on $[-n^{\tau}, n^{\tau}]$ for some fixed $0<\tau<1/(d+1)$.
Then we set $\alpha_n = c_1 \cdot \log n$, and compute
\[
\bw^{(t+1)}=\bw^{(t)} - \lambda_n \cdot (\nabla_\bw F_n)(\bw^{(t)})
\quad (t=0, \dots, t_n-1)
\]
where
\[
F_n(\bw) = \frac{1}{n} \sum_{i=1}^n |f_\bw(X_i) -Y_i|^2 \cdot 1_{[-\alpha_n,\alpha_n]^d}(X_i)
+ c_2 \cdot \sum_{j=1}^{K_n} |w_{1,1,j}^{(L)}|^2
\]
is the regularized empirical $L_2$ risk of the network $f_\bw$ on the training data.
The step size $\lambda_n>0$ and the number $t_n$
of gradient descent steps will be chosen in Theorem \ref{th1}
below.

The estimate is defined by
\[
m_n(x) = (T_{\beta_n} f_{\bw^{(t_n)}} (x)) \cdot 1_{[-\alpha_n,\alpha_n]^d}(x),
\]
i.e., as an estimate we use the neural network with the weight
vector which we get after $t_n$ gradient descent steps, and truncate
this function on height $-\beta_n$ and $\beta_n$, and set it
equal to zero outside of a cube. Here
we set $\beta_n = c_3 \cdot \log n$.

Because of (\ref{se2eq5}) we have
\[
F_n(\bw^{(0)}) = \frac{1}{n} \sum_{i=1}^n |Y_i|^2
\cdot 1_{[-\alpha_n,\alpha_n]^d}(X_i).
\]

\section{Main result}
\label{se3}

Our main result is the following theorem.

\begin{theorem}
  \label{th1}
  Let $\sigma$ be the logistic squasher, and let
  $K_n, L, r \in \N$ and $\tau \in \R_+$.
  Assume $L \geq 2$, $r \geq 2d$, $0<\tau<1/(d+1)$
  \begin{equation}
    \label{th1eq1}
    \frac{K_n}{n^\kappa} \rightarrow 0 \quad (n \rightarrow \infty)
  \end{equation}
  for some $\kappa>0$,
  \begin{equation}
    \label{th1eq2}
   \frac{K_n}{n^{r} \cdot \log n} \rightarrow \infty \quad (n \rightarrow \infty), 
  \end{equation}
  and set $\alpha_n= c_1 \cdot \log n$, $\beta_n= c_3 \cdot \log n$,
  \begin{equation}
    \label{th1eq3}
    \lambda_n= \frac{1}{L_n} \quad \mbox{and} \quad
    t_n = \lceil c_4 \cdot L_n \cdot \log n \rceil
  \end{equation}
  for some $L_n>0$ which satisfies
  \begin{equation}
    \label{th1eq4}
    L_n \geq (\log n)^{10 \cdot L + 10} \cdot K_n^{3/2}.
  \end{equation}
  Let the estimate $m_n$ be defined as in Section \ref{se2}.
  Then we have
  \[
  \EXP \int | m_n(x)-m(x)|^2 \PROB_X (dx)
  \rightarrow 0
  \quad
  (n \rightarrow \infty)
  \]
  for every distribution of $(X,Y)$ with $\EXP Y^2 < \infty$.
  \end{theorem}

\noindent
    {\bf Remark 1.} Condition (\ref{th1eq2}) implies that
    $K_n$ is asymptotically larger than $n^{r}$, consequently
    the number of parameters of
    our estimate is much larger than the sample size
    and our estimate is over-parametrized. Nevertheless it
    generalizes well on new independent data since its expected
    $L_2$ error converges to zero for sample size tending
    to infinity. In its definition we add a regularization term
    to the empirical $L_2$ risk,
    however, this regularization term is not really used
    to control the complexity of our estimate, it is only used to
    help us in analyzing the gradient descent. We control the
    complexity of our estimate by imposing bounds on the absolute values
    of the initial weights and by requiring that the number of gradient
    descent steps is not much larger than the reciprocal of the stepsize.
    In particular, the condition $0 < \tau < 1/(d+1)$ controls the range
    $[-n^\tau,n^\tau]$ of the weights $w_{k,i,j}^{(0)}$, and all initial weights
    of level $1, \dots, L-1$ are bounded in absolute value
    by some logarithmic term.
    
\noindent
    {\bf Remark 2.} We need only a single initialization of our
    random starting weights in Theorem \ref{th1}. This is due to the
    over-parametrization, which enables us to show that even with
    one single initialization there exists with probability close to
    one a finite subset of our $K_n$ fully connected neural networks
    where the initial inner weights have some nice property.

    \noindent
        {\bf Remark 3.} In the proof of Theorem \ref{th1} we use that
        the inner weights do not change much during gradient descent
        and that gradient descent is able to find proper values for
        the outer weights of our network. In this sense our network
        is not based on representation learning, instead it is
        using a representation guessing.

\section{Proofs}
\label{se4}

    \subsection{Auxiliary results}
    \label{se4sub1}

    In this subsection we present various auxiliary results which we
    will need in the proof of Theorem \ref{th1}.
    
\begin{lemma}
  \label{le1}
  Let
  $F:\R^K \rightarrow \R_+$
  be a nonnegative differentiable function.
  Let
  $t \in \N$, $L>0$, $\ba_0 \in \R^K$ and set
  \[
  \lambda=
\frac{1}{L}
\]
and
\[
\ba_{k+1}=\ba_k - \lambda \cdot (\nabla_{\ba} F)(\ba_k)
\quad
(k \in \{0,1, \dots, t-1\}).
\]
Assume
\begin{equation}
  \label{le1eq1}
  \left\|
 (\nabla_{\ba} F)(\ba)
  \right\|
  \leq
  \sqrt{
2 \cdot t \cdot L \cdot \max\{ F(\ba_0),1 \}
    }
\end{equation}
for all $\ba \in \R^K$ with
$\| \ba - \ba_0\| \leq \sqrt{2 \cdot t \cdot \max\{ F(\ba_0),1 \} / L}$,
and
\begin{equation}
  \label{le1eq2}
\left\|
(\nabla_{\ba} F)(\ba)
-
(\nabla_{\ba} F)(\bb)
  \right\|
  \leq
  L \cdot \|\ba - \bb \|
\end{equation}
for all $\ba, \bb \in \R^K$ satisfying
\begin{equation}
  \label{le1eq3}
  \| \ba - \ba_0\| \leq \sqrt{8 \cdot \frac{t}{L} \cdot \max\{ F(\ba_0),1 \}}
  \quad \mbox{and} \quad
  \| \bb - \ba_0\| \leq \sqrt{8 \cdot \frac{t}{L} \cdot \max\{ F(\ba_0),1 \}}.
\end{equation}
Then we have
\[
\|\ba_k-\ba_0\| \leq
\sqrt{
2 \cdot \frac{k}{L} \cdot (F(\ba_0)-F(\ba_k))
}
\quad
 \mbox{for all }
 k \in \{1, \dots,t\},
\]
\[
\sum_{k=0}^{s-1}
\| \ba_{k+1}-\ba_k \|^2
\leq
\frac{2}{L}
 \cdot (F(\ba_0)-F(\ba_s))
\quad
 \mbox{for all }
 s \in \{1, \dots,t\}
 \]
 and
 \[
 F(\ba_k)
 \leq
 F(\ba_{k-1})
                -
                \frac{1}{2 L} \cdot
                \| \nabla_{\ba}  F(\ba_{k-1}) \|^2
                \quad
 \mbox{for all }
 k \in \{1, \dots,t\}.
\]
\end{lemma}

\noindent
    {\bf Proof.} The result follows from Lemma 2 in Braun et al. (2021)
    and its proof.
    \hfill $\Box$

\begin{lemma}\label{le2}
  Let $\sigma: \R \rightarrow \R$ be bounded and differentiable, and assume that
its derivative is bounded.
Let $\alpha_n \geq 1$,
$t_n \geq L_n$,
$\gamma_n^* \geq 1$, $B_n \geq 1$, $r \geq 2d$,
\begin{equation}
	\label{le2eq1}
	|w_{1,1,k}^{(L)}| \leq \gamma_n^* \quad (k=1, \dots, K_n),
\end{equation}
\begin{equation}
	\label{le2eq2}
	|w_{k,i,j}^{(l)}| \leq B_n
	\quad
	\mbox{for } l=1, \dots, L-1
\end{equation}
and
\begin{equation}
	\label{le2eq3}
	\|\bw-\bv\|_\infty^2 \leq \frac{2t_n}{L_n} \cdot \max\{ F_n(\bv),1 \}.
\end{equation}
Then we have 
\[
\| (\nabla_\bw F_n)(\bw) \|
\leq
c_5 \cdot K_n^{3/2} \cdot B_n^{2L} \cdot (\gamma_n^*)^2 \cdot \alpha_n^{2} \cdot \sqrt{\frac{t_n}{L_n} \cdot \max\{F_n(\bv),1\}}. 
\]
\end{lemma}

\noindent
    {\bf Proof.}
    We have
    \begin{eqnarray*}
      &&
      \| \nabla_\bw F_n (\bw) \|^2
      \\
      &&
      =
      \sum_{k,i,j,l}
      \Bigg(
      \frac{2}{n}
      \sum_{s=1}^n
      (Y_s - f_\bw (X_s)) \cdot 1_{[-\alpha_n,\alpha_n]^d}(X_s)
      \cdot
      \frac{\partial f_\bw}{\partial w_{k,i,j}^{(l)}}(X_s)
      \\
      &&
      \hspace*{3cm}
      +
      \frac{\partial}{\partial w_{k,i,j}^{(l)}}
      \left(
      c_2 \cdot
      \sum_{r=1}^{K_n} |w_{1,1,r}^{(L)}|^2
      \right)
      \Bigg)^2
      \\
      &&
      \leq
      8 \cdot
      \sum_{k,i,j,l}
      \frac{1}{n}
      \sum_{s=1}^n
      (Y_s - f_\bw (X_s))^2 \cdot 1_{[-\alpha_n,\alpha_n]^d}(X_s)
      \cdot
      \left(
      \frac{\partial f_\bw}{\partial w_{k,i,j}^{(l)}}(X_s)
      \right)^2
      \\
      &&
      \hspace*{3cm}
      + 8 \cdot c_2^2 \cdot K_n \cdot (\gamma_n^*)^2
      \\
      &&
      \leq
      c_6 \cdot K_n \cdot L \cdot r^2 \cdot d \cdot
      \max_{k,i,j,l,s}       \left(
      \frac{\partial f_\bw}{\partial w_{k,i,j}^{(l)}}(X_s)
      \right)^2
      \cdot 1_{[-\alpha_n,\alpha_n]^d}(X_s)
      \\
      &&
      \hspace*{3cm}
      \cdot
      \frac{1}{n}
      \sum_{s=1}^n
      (Y_s - f_\bw (X_s))^2 \cdot 1_{[-\alpha_n,\alpha_n]^d}(X_s)
      + 8 \cdot c_2^2 \cdot K_n \cdot (\gamma_n^*)^2.
    \end{eqnarray*}
    The chain rule implies
    \begin{eqnarray}
      &&
      \frac{\partial f_\bw}{\partial w_{k,i,j}^{(l)}}(x)
        =
  \sum_{s_{l+2}=1}^{r} \dots \sum_{s_{L-1}=1}^{r}
  f_{k,j}^{(l)}(x)
  \cdot
\sigma^\prime \left(\sum_{t=1}^{r} w_{k,i,t}^{(l)} \cdot f_{k,t}^{(l)}(x) + w_{k,i,0}^{(l)} \right)
  \nonumber \\
  && \quad
  \cdot
  w_{k,s_{l+2},i}^{(l+1)} \cdot
\sigma^\prime \left(\sum_{t=1}^{r} w_{k,s_{l+2},t}^{(l+1)} \cdot f_{k,t}^{(l+1)}(x) + w_{k,s_{l+2},0}^{(l+1)} \right)
  \cdot
  w_{k,s_{l+3},s_{l+2}}^{(l+2)}
   \nonumber \\
  && \quad
  \cdot
\sigma^\prime \left(\sum_{t=1}^{r} w_{k,s_{l+3},t}^{(l+2)} \cdot f_{k,t}^{(l+2)}(x) + w_{k,s_{l+3},0}^{(l+2)} \right)
  \cdots
  w_{k,s_{L-1},s_{L-2}}^{(L-2)}
 \nonumber \\
  && \quad
  \cdot
\sigma^\prime \left(\sum_{t=1}^{r} w_{k,s_{L-1},t}^{(L-2)} \cdot f_{k,t}^{(L-2)}(x) + w_{k,s_{L-1},0}^{(L-2)} \right)
   \cdot
   w_{k,1,s_{L-1}}^{(L-1)}
 \nonumber \\
  && \quad
   \cdot
\sigma^\prime \left(\sum_{t=1}^{r} w_{k,1,t}^{(L-1)} \cdot f_{k,t}^{(L-1)}(x) + w_{k,1,0}^{(L-1)} \right)
     \cdot
     w_{1,1,k}^{(L)},
     \label{ple2eq1}
      \end{eqnarray}
where we have used the abbreviations
\[
f_{k,j}^{(0)}(x)
=
\left\{
\begin{array}{ll}
  x^{(j)} & \mbox{if } j \in \{1,\dots,d\} \\
  1 & \mbox{if } j=0
\end{array}
\right.
\]
and
\[
f_{k,0}^{(l)}(x)=1 \quad (l=1, \dots, L-1).
\]
Using the assumptions of Lemma \ref{le2} we can conclude
\[
      \max_{k,i,j,l,s}       \left(
      \frac{\partial f_\bw}{\partial w_{k,i,j}^{(l)}}(X_s)
      \right)^2
\cdot 1_{[-\alpha_n,\alpha_n]^d}(X_s)      
      \leq
      c_6 \cdot r^{2L} \cdot
      \max\{ \|\sigma^\prime\|_\infty^{2L},1\} \cdot B_n^{2L} \cdot (\gamma_n^*)^2 \cdot \alpha_n^2.
      \]
      By the Lipschitz continuity of $\sigma$ together with the
      assumptions of Lemma \ref{le2} we get for any $x \in [-\alpha_n,\alpha_n]^d$
      \[
      |f_\bw(x)-f_\bv(x)| \leq
     2\cdot K_n \cdot \max\{\|\sigma'\|_\infty^L,1\} \cdot \gamma_n^* \cdot (2r+1)^L \cdot B_n^L \cdot \alpha_n \cdot \max\{\|\sigma\|_\infty,1\} \cdot \|\bw-\bv\|_\infty.
      \]
(cf., e.g., Lemma 5 in Kohler and Krzy\.zak (2021) for a related proof).
      This implies
      \begin{eqnarray*}
        &&
      \frac{1}{n}
      \sum_{s=1}^n
      (Y_s - f_\bw (X_s))^2 \cdot 1_{[-\alpha_n,\alpha_n]^d}(X_s)
      \\
      &&
      \leq
      2 \cdot F_n(\bv) +       \frac{2}{n}
      \sum_{s=1}^n
      (f_\bv(X_s) - f_\bw (X_s))^2 \cdot 1_{[-\alpha_n,\alpha_n]^d}(X_s)
      \\
      &&
      \leq 2 \cdot F_n(\bv) + 8 \cdot \max\{\|\sigma'\|_\infty^{2L},1\} \cdot K_n^2 \cdot {\gamma_n^*}^2 \cdot (2r+1)^{2L} \cdot B_n^{2L} \cdot \alpha_n^{2} \cdot \max\{\|\sigma\|_\infty,1\}^2
      \\
      &&
      \quad
      \cdot \frac{2t_n}{L_n} \cdot \max\{F_n(\bv),1\}.
        \end{eqnarray*}
Summarizing the above results, the proof is complete.       
    \hfill $\Box$

    \begin{lemma}
      \label{le3}
     Let $\sigma: \R \rightarrow \R$ be bounded and differentiable, and assume that its derivative
     is 
     Lipschitz continuous and bounded.
     Let $\alpha_n \geq 1$, 
     $t_n \geq L_n$,
     $\gamma_n^* \geq 1$, $B_n \geq 1$, $r \geq 2d$ and assume
     \begin{equation}
     	\label{le3eq1}
     	|\max\{ (\bw_1)_{1,1,k}^{(L)}, (\bw_2)_{1,1,k}^{(L)}\}| \leq \gamma_n^* \quad (k=1,
     	\dots, K_n),
     \end{equation}
     
     \begin{equation}
     	\label{le3eq2}
     	|\max\{(\bw_1)_{k,i,j}^{(l)},(\bw_2)_{k,i,j}^{(l)}\}| \leq B_n
     	\quad
     	\mbox{for } l=1, \dots, L-1
     \end{equation}
     and
     \begin{equation}
     	\label{le3eq3}
     	\|\bw_2-\bv\|^2 \leq 8 \cdot \frac{t_n}{L_n} \cdot \max\{ F_n(\bv),1 \}.
     \end{equation}
     Then we have 
     \begin{eqnarray*}
     	&&
     	\| (\nabla_\bw F_n)(\bw_1) - (\nabla_\bw F_n)(\bw_2) \| \\
     	&&
     	\leq
     	c_7 \cdot \max \{\sqrt{F_n(\bv)},1\} \cdot (\gamma_n^*)^{2} \cdot B_n^{3L} \cdot \alpha_n^{3} \cdot K_n^{3/2} \cdot \sqrt{\frac{t_n}{L_n}} \cdot \|\bw_1-\bw_2\|. 
     \end{eqnarray*}
\end{lemma}

\noindent
    {\bf Proof.}
We have
    \begin{eqnarray*}
      &&
      \| \nabla_\bw F_n (\bw_1) -  \nabla_\bw F_n (\bw_2) \|^2
      \\
      &&
      =
      \sum_{k,i,j,l}
      \Bigg(
      \frac{2}{n}
      \sum_{s=1}^n
      (Y_s - f_{\bw_1} (X_s)) \cdot 1_{[-\alpha_n,\alpha_n]^d}(X_s)
      \cdot
      \frac{\partial f_{\bw_1}}{\partial w_{k,i,j}^{(l)}}(X_s)
      \\
      &&
      \hspace*{3cm}
      +
      \frac{\partial}{\partial w_{k,i,j}^{(l)}}
      \left(
      c_2 \cdot
      \sum_{r=1}^{K_n} |(\bw_1)_{1,1,r}^{(L)}|^2
      \right)
      \\
      &&
      \quad
      -
      \sum_{k,i,j,l}
      \Bigg(
      \frac{2}{n}
      \sum_{s=1}^n
      (Y_s - f_{\bw_2} (X_s)) \cdot 1_{[-\alpha_n,\alpha_n]^d}(X_s)
      \cdot
      \frac{\partial f_{\bw_2}}{\partial w_{k,i,j}^{(l)}}(X_s)
      \\
      &&
      \hspace*{3cm}
      +
      \frac{\partial}{\partial w_{k,i,j}^{(l)}}
      \left(
      c_2 \cdot
      \sum_{r=1}^{K_n} |(\bw_2)_{1,1,r}^{(L)}|^2
      \right)      
      \Bigg)^2
      \\
      &&
      \leq
      16 \cdot
      \sum_{k,i,j,l}
      \left(
      \frac{1}{n}
      \sum_{s=1}^n
      (f_{\bw_2} (X_s) - f_{\bw_1} (X_s)) \cdot 1_{[-\alpha_n,\alpha_n]^d}(X_s)
      \cdot
      \frac{\partial f_{\bw_1}}{\partial w_{k,i,j}^{(l)}}(X_s)
      \right)^2
      \\
      &&
      \quad
      +
      16 \cdot
            \sum_{k,i,j,l}
      \left(
      \frac{1}{n}
      \sum_{s=1}^n
      (Y_s-f_{\bw_2} (X_s) ) \cdot 1_{[-\alpha_n,\alpha_n]^d}(X_s)
      \cdot
      \left(
      \frac{\partial f_{\bw_1}}{\partial w_{k,i,j}^{(l)}}(X_s)
      -
           \frac{\partial f_{\bw_2}}{\partial w_{k,i,j}^{(l)}}(X_s)
 \right)
      \right)^2
      \\
      &&
      \quad
      +
      8 \cdot c_2^2 \cdot \|\bw_1-\bw_2\|^2
      \\
      &&
      \leq
      16 \cdot
      \sum_{k,i,j,l}
      \max_{s=1,\dots,n}
      \left(
      \frac{\partial f_{\bw_1}}{\partial w_{k,i,j}^{(l)}}(X_s)
      \right)^2 \cdot 1_{[-\alpha_n,\alpha_n]^d}(X_s)
  \\
      &&
      \hspace*{3cm}
    \cdot
      \frac{1}{n}
      \sum_{s=1}^n
      (f_{\bw_2} (X_s) - f_{\bw_1} (X_s))^2
      \cdot 1_{[-\alpha_n,\alpha_n]^d}(X_s)
      \\
      &&
      \quad
      +
      16  \cdot
      \frac{1}{n}
      \sum_{s=1}^n
      (Y_s-f_{\bw_2} (X_s) )^2 \cdot 1_{[-\alpha_n,\alpha_n]^d}(X_s)
      \\
      &&
      \hspace*{3cm}
      \cdot
      \sum_{k,i,j,l}
      \max_{s=1,\dots,n}
      \left(
      \frac{\partial f_{\bw_1}}{\partial w_{k,i,j}^{(l)}}(X_s)
      -
           \frac{\partial f_{\bw_2}}{\partial w_{k,i,j}^{(l)}}(X_s)
           \right)^2 \cdot 1_{[-\alpha_n,\alpha_n]^d}(X_s)
           \\
           &&
           \quad
      +
      8 \cdot c_2^2 \cdot \|\bw_1-\bw_2\|_\infty^2
      .
    \end{eqnarray*}
    From the proof of Lemma \ref{le2} we can conclude
    \begin{eqnarray*}
      &&
      \sum_{k,i,j,l}
      \max_{s=1,\dots,n}
      \left(
      \frac{\partial f_{\bw_1}}{\partial w_{k,i,j}^{(l)}}(X_s)
      \right)^2 \cdot 1_{[-\alpha_n,\alpha_n]^d}(X_s)
      \\
      &&
      \leq
      c_8 \cdot K_n \cdot L \cdot r^2 \cdot d \cdot
      r^{2L} \cdot
      \max\{ \|\sigma^\prime\|_\infty^{2L},1 \} \cdot B_n^{2L} \cdot (\gamma_n^*)^2 \cdot \alpha_n^2,      
      \end{eqnarray*}
    \begin{eqnarray*}
      &&
      \frac{1}{n}
      \sum_{s=1}^n
      (f_{\bw_2} (X_s) - f_{\bw_1} (X_s))^2 \cdot 1_{[-\alpha_n,\alpha_n]^d}(X_s)
      \\
      &&
      \leq
      4 \cdot \max\{\|\sigma'\|_\infty^{2L},1\} \cdot K_n^2 \cdot (2r+1)^{2L} \cdot
      (\gamma_n^*)^2 
      \cdot B_n^{2L} \cdot \alpha_n^2 \cdot \max\{\|\sigma\|_\infty,1\}^2 \cdot \|\bw_1-\bw_2\|^2
      \end{eqnarray*}
    and
    \begin{eqnarray*}
      &&
            \frac{1}{n}
      \sum_{s=1}^n
      (Y_s-f_{\bw_2} (X_s) )^2 \cdot 1_{[-\alpha_n,\alpha_n]^d}(X_s)
      \\
      &&
      \leq
            2 \cdot F_n(\bv)
      + 8 \cdot
      \max\{\|\sigma'\|_\infty^{2L},1\} \cdot K_n^2 \cdot (2r+1)^{2L} \cdot (\gamma_n^*)^2
      \cdot B_n^{2L} \cdot \alpha_n^2 \cdot \max\{\|\sigma\|_\infty,1\}^2
      \\
      &&
      \quad
      \cdot
      \frac{8 t_n}{L_n} \cdot \max\{F_n(v),1\}.
    \end{eqnarray*}
    So it remains to bound
    \[
\sum_{k,i,j,l}
      \max_{s=1,\dots,n}
      \left(
      \frac{\partial f_{\bw_1}}{\partial w_{k,i,j}^{(l)}}(X_s)
      -
           \frac{\partial f_{\bw_2}}{\partial w_{k,i,j}^{(l)}}(X_s)
 \right)^2 \cdot 1_{[-\alpha_n,\alpha_n]^d}(X_s).
 \]
 By (\ref{ple2eq1}) we know that
 \[
 \frac{\partial f_{\bw}}{\partial w_{k,i,j}^{(l)}}(x)
 \]
 is for fixed $x \in [-\alpha_n,\alpha_n]^d$ a sum of products
 of Lipschitz continuous functions (considered as functions of $\bw$).
 Arguing as in the proof of Lemma 6 in Kohler and Krzy\.zak (2021)
 we can show that we have for any $x \in [-\alpha_n,\alpha_n]^d$
 \begin{eqnarray*}
   &&
   \left|
   \frac{\partial f_{\bw_1}}{\partial w_{k,i,j}^{(l)}}(x)
      -
           \frac{\partial f_{\bw_2}}{\partial w_{k,i,j}^{(l)}}(x)
           \right|
           \leq
           c_9 \cdot B_n^{2L} \cdot \gamma_n^* \cdot \alpha_n \cdot \|\bw_1-\bw_2\|,
   \end{eqnarray*}
 which implies
 \begin{eqnarray*}
   &&
\sum_{k,i,j,l}
      \max_{s=1,\dots,n}
      \left(
      \frac{\partial f_{\bw_1}}{\partial w_{k,i,j}^{(l)}}(X_s)
      -
           \frac{\partial f_{\bw_2}}{\partial w_{k,i,j}^{(l)}}(X_s)
           \right)^2 \cdot 1_{[-\alpha_n,\alpha_n]^d}(X_s)
           \\
           &&
           \leq
      c_{10} \cdot K_n \cdot L \cdot r^2 \cdot d \cdot
  B_n^{4L} \cdot (\gamma_n^*)^2 \cdot \alpha_n^4 \cdot \|\bw_1-\bw_2\|^2.         
 \end{eqnarray*}
 Summarizing the above results we get the assertion.
    \hfill $\Box$

\begin{lemma}
  \label{le4}
  Let $\alpha \geq 1$, $\beta>0$ and let $A,B,C \geq 1$.
  Let $\sigma:\R \rightarrow \R$ be $k$-times differentiable
  such that all derivatives up to order $k$ are bounded on $\R$.
  Let $\F$
  be the set of all functions $f_{\bw}$ defined by
  (\ref{se2eq1})--(\ref{se2eq3}) where the weight vector $\bw$
  satsifies
  \begin{equation}
    \label{le5eq1}
    \sum_{j=1}^{K_n} |w_{1,1,j}^{(L)}| \leq C,
    \end{equation}
  \begin{equation}
    \label{le4eq2}
    |w_{k,i,j}^{(l)}| \leq B \quad (k \in \{1, \dots, K_n\},
    i,j \in \{1, \dots, r\}, l \in \{1, \dots, L-1\})
    \end{equation}
and
  \begin{equation}
    \label{le4eq3}
    |w_{k,i,j}^{(0)}| \leq A \quad (k \in \{1, \dots, K_n\},
    i \in \{1, \dots, r\}, j \in \{1, \dots,d\}).
  \end{equation}
  Then we have for any $1 \leq p < \infty$, $0 < \epsilon < \beta$ and
  $x_1^n \in \Rd$
  \begin{eqnarray*}
&&\Nu_p \left(
\epsilon, \{ T_\beta f \cdot 1_{[-\alpha,\alpha]^d} \, : \, f \in \F \}, x_1^n
\right)
\\
&&
\leq \left(c_{11}\cdot \frac{\beta^p} {\epsilon^p}\right)^{c_{12}\cdot \alpha^d \cdot B^{(L-1)\cdot d} \cdot A^d \cdot \left(\frac{C}{\epsilon}\right)^{d/k}+ c_{13}
}.
\\
  \end{eqnarray*}
  
  \end{lemma}

\noindent
    {\bf Proof.} In the {\it first step} of the proof we show
    for any $f_{\bw} \in \F$, any $x \in \Rd$ and any $s_1, \dots, s_k \in \{1, \dots, d\}$
    \begin{equation}
      \label{ple4eq1}
      \left|
      \frac{\partial^k f_{\bw}}{\partial x^{(s_1)} \dots \partial x^{(s_k)}} (x)
      \right| \leq c_{14} \cdot C \cdot B^{(L-1) \cdot k} \cdot A^k =: c
      .
      \end{equation}
    The definition of $f_{\bw}$ implies
    \[
    \frac{\partial^k f_{\bw}}{\partial x^{(s_1)} \dots \partial x^{(s_k)}} (x) = \sum_{j=1}^{K_n} w_{1,1,j}^{(L)} \cdot
    \frac{\partial^k
f_{j,1}^{(L)}(x)
    }{\partial x^{(s_1)} \dots \partial x^{(s_k)}} (x),
    \]
    hence (\ref{ple4eq1}) is implied by
    \begin{equation}
      \label{ple4eq1b}
      \left|
      \frac{\partial^k f_{j,1}^{(L)}}{\partial x^{(s_1)} \dots \partial x^{(s_k)}} (x)
      \right| \leq c_{15} \cdot B^{(L-1) \cdot k} \cdot A^k
      .
      \end{equation}
We have
    \begin{eqnarray*}
      \frac{\partial f_{k,i}^{(l)}}{\partial x^{(s)}}(x)
&
      =
&
      \sigma^\prime \left(\sum_{t=1}^{r} w_{k,i,t}^{(l-1)} \cdot f_{k,t}^{(l-1)}(x) + w_{k,i,0}^{(l-1)} \right)
      \cdot
      \sum_{j=1}^{r} w_{k,i,j}^{(l-1)} \cdot \frac{\partial
f_{k,j}^{(l-1)}
      }{\partial x^{(s)}}(x)
      \\
      &
      =&
      \sum_{j=1}^{r} w_{k,i,j}^{(l-1)} 
      \cdot
      \sigma^\prime \left(\sum_{t=1}^{r} w_{k,i,t}^{(l-1)} \cdot f_{k,t}^{(l-1)}(x) + w_{k,i,0}^{(l-1)} \right)
      \cdot
      \frac{\partial
f_{k,j}^{(l-1)}
      }{\partial x^{(s)}}(x)
    \end{eqnarray*}
    and
    \begin{eqnarray*}
      \frac{\partial f_{k,i}^{(1)}}{\partial x^{(s)}}(x)
      &=&
      \sigma^\prime \left(\sum_{j=1}^d w_{k,i,j}^{(0)} \cdot x^{(j)} + w_{k,i,0}^{(0)} \right)
      \cdot
       w_{k,i,s}^{(0)}.
      \end{eqnarray*}
    By the product rule of derivation we can conclude for $l>1$ that
    \begin{equation}
      \label{ple4eq2}
    \frac{\partial^k f_{k,i}^{(l)}}{\partial x^{(s_1)} \dots \partial x^{(s_k)}} (x)
    \end{equation}
    is a sum of at most $r \cdot (r+k)^{k-1}$ terms of the form
    \begin{eqnarray*}
      &&
w
\cdot
\sigma^{(s)} \left(\sum_{j=1}^{r} w_{k,i,j}^{(l-1)} \cdot f_{k,j}^{(l-1)}(x) + w_{k,i,0}^{(l-1)} \right)
\\
&&
\hspace*{3cm}
      \cdot
      \frac{\partial^{t_{1}}
f_{k,j_1}^{(l-1)}
      }{\partial x^{(r_{1,1})} \dots \partial x^{(r_{1,t_{1}})}}(x)
      \cdot
\dots
      \cdot
      \frac{\partial^{t_{s}}
f_{k,j_s}^{(l-1)}
      }{\partial x^{(r_{s,1})} \dots \partial x^{(r_{s,t_{s}})}}(x)
    \end{eqnarray*}
    where we have $s \in \{1, \dots,k\}$,
    $|w| \leq  B^{s}$ and $t_1+ \dots + t_s = k$.
    Furthermore
    \[
    \frac{\partial^k f_{k,i}^{(1)}}{\partial x^{(s_1)} \dots \partial x^{(s_k)}} (x)
    \]
    is a given by
    \[
    \prod_{j=1}^k w_{k,i,s_j}^{(0)} \cdot
      \sigma^{(k)} \left(\sum_{t=1}^d w_{k,i,t}^{(0)} \cdot x^{(t)} + w_{k,i,0}^{(0)} \right)
.
    \]
    Because of the boundedness of the derivatives of $\sigma$ we can conclude
    from (\ref{le4eq3})
    \[
    \left|
    \frac{\partial^k f_{k,i}^{(1)}}{\partial x^{(s_1)} \dots \partial x^{(s_k)}} (x)
\right|
\leq
c_{16} \cdot
A^k
    \]
    for all $k \in \N$ and $s_1, \dots, s_k \in \{1, \dots, d\}$.

    Recursively we can conclude from the above
    representation of (\ref{ple4eq2}) that we have
    \[
    \left|
    \frac{\partial^k f_{k,i}^{(l)}}{\partial x^{(s_1)} \dots \partial x^{(s_k)}} (x)
    \right|
    \leq
    c_{17} \cdot B^{(l-1) \cdot k} \cdot A^k.
    \]
    Setting $l=L$ we get (\ref{ple4eq1b}).

In the {\it second step} of the proof we show
\begin{equation}
	\label{ple4eq3}
	\Nu_p \left(
	\epsilon, \{T_\beta f \cdot 1_{[-\alpha,\alpha]^d} \, : \, f \in \F \}, x_1^n
	\right)
	\leq
	\Nu_p  \left(
	\frac{\epsilon}{2}, T_\beta \G\circ \Pi, x_1^n
	\right),     
\end{equation}
where $\G$ is the set of all polynomials of degree less than or equal to $k-1$
which vanish outside of $[-\alpha,\alpha]^d$
and $\Pi$ is the family of all partitions of $\mathbb{R}^d$ which consist of a partition of $[-\alpha,\alpha]^d$ into $K$
many cubes of sidelenght 
\[
\left(c_{18}\cdot \frac{\epsilon}{c}\right)^{1/k}
\]
where $c_{18}=c_{18}(d,k)$ is a suitable small constant greater than zero
and the additional set $\mathbb{R}^d\setminus [-\alpha,\alpha]^d$.

A standard bound on the remainder of a multivariate Taylor polynomial
together with (\ref{ple4eq1})
shows that for each $f_\bw$ we can find $g \in \mathcal{G}\circ \Pi$ such that
\[
|f_\bw(x)-g(x)|\leq \frac{\epsilon}{2}
\]
  holds for all $x\in [-\alpha,\alpha]^d$, which implies (\ref{ple4eq3}).

    In the {\it third step} of the proof we show the assertion of Lemma \ref{le4}.
    Since $\mathcal{G}\circ \Pi$ is a linear vector space of dimension less than or equal to \[
    c_{20} \cdot \alpha^d \cdot \left(\frac{c}{\epsilon}\right)^{d/k}
    \]
    we conclude from Theorem 9.4 and Theorem 9.5 in Gy\"orfi et al. (2002),
    	\begin{align*}
    	\mathcal{N}_p(\frac{\epsilon}{2}, T_\beta \mathcal{G} \circ \Pi, x_1^n)
    	\leq 3 \left(\frac{2e(2\beta)^p}{(\epsilon/2)^p}\log\left(\frac{3e(2 \beta)^p}{(\epsilon/2)^p}\right)\right)^{c_{20} \cdot \alpha^d \cdot \left(\frac{c}{\epsilon}\right)^{d/k}+1}.
    \end{align*}

    Together with (\ref{ple4eq3}) this implies the assertion.
    \quad  \hfill $\Box$

\begin{lemma}
  \label{le5}
Let $\sigma$ be the logistic squasher and let $0<\delta \leq 1$, $1 \leq \alpha_n \leq \log n$, $\mathbf{u},\mathbf{v} \in \mathbb{R}^d$ with
\begin{equation*}
v^{(l)}-u^{(l)} \geq 2\delta \quad \mbox{for } l\in\{1,\dots, d\}
\end{equation*}
and $x \in
  [-\alpha_n,\alpha_n]^d$.
Let $L, r,n \in \mathbb{N}$ with $L \geq 2$, $r \geq 2 \cdot d$, $n \geq 8d$ and $n\geq \exp(r+1)$.
Let
\[
f_{\bw}(x)= f_{1,1}^{(L)}(x) 
\]
where $f_{k,i}^{(l)}(x)$ are recursively defined by (\ref{se2eq2}) and (\ref{se2eq3}).

Assume
\begin{equation}
\label{le5eq1}
w_{1,j,j}^{(0)}=\frac{4d\cdot (\log n)^2}{\delta}
\quad \mbox{and} \quad  w_{1,j,0}^{(0)}=-\frac{4d \cdot (\log n)^2}{\delta} \cdot u^{(j)} \quad \text{for }j \in \{1, \dots, d\},
\end{equation}
\begin{equation}
\label{le5eq2}
w_{1,j+d,j}^{(0)}=-\frac{4d \cdot (\log n)^2}{\delta}
 \quad \mbox{and}\quad w_{1,j+d,0}^{(0)}=\frac{4d \cdot (\log n)^2}{\delta} \cdot v^{(j)}
\quad \text{for }j \in \{1, \dots, d\},
\end{equation}
\begin{equation}
\label{le5eq3}
w_{1,s,t}^{(0)}=0
\quad \mbox{if }
s \leq 2d,
s \neq t, s \neq t+d \mbox{ and } t>0,
\end{equation}
\begin{equation}
\label{le5eq4}
w_{1,1,t}^{(1)}= 8 \cdot (\log n)^2 \quad \mbox{for }t \in \{1, \dots, 2d\},
\end{equation}
\begin{equation}
\label{le5eq5}
w_{1,1,0}^{(1)} = - 8 \cdot (\log n)^2\left(2d-\frac{1}{2}\right),
\end{equation}
\begin{equation}
\label{le5eq6}
w_{1,1,t}^{(1)}=0 
\quad \mbox{for }t > 2d,
\end{equation}
\begin{equation}
\label{le5eq7}
w_{1,1,1}^{(l)}= 6 \cdot (\log n)^2 \quad \mbox{for }l \in \{2, \dots, L\},
\end{equation}
\begin{equation}
\label{le5eq8}
w_{1,1,0}^{(l)} = - 3\cdot(\log n)^2 \quad \mbox{for }l \in \{2, \dots, L\}
\end{equation}
and
\begin{equation}
\label{le5eq9}
w_{1,1,t}^{(l)}=0 
\quad \mbox{for }
 t > 1\mbox{ and } l \in \{2, \dots, L\}.
\end{equation}

Let $\bar{\bw}$ be such that 
\begin{equation}
\label{le5eq10}
|\bar{w}_{1,i,j}^{(l)}-w_{1,i,j}^{(l)}| \leq \log n
\quad \mbox{for all } l=0,\dots,L-1.
\end{equation}
Then, we have
\[
f_{\bar{\bw}}(x) \geq 1 - \frac{1}{n} \mbox{ if } x \in
[u^{(1)}+\delta, v^{(1)}-\delta] \times \dots \times
[u^{(d)}+\delta, v^{(d)}-\delta] 
\]
and
\[
f_{\bar{\bw}}(x) \leq \frac{1}{n} \mbox{ if } x^{(i)} \notin
[u^{(i)}-\delta, v^{(i)}+\delta] \mbox{ for some } i \in \{1,\dots, d\}. 
\]
  \end{lemma}

\noindent
    {\bf Proof.}
At the beginning we define

\begin{equation*}
\bar{f}_{k,i}^{(l)}(x) = \sigma\left(\sum_{j=1}^{r} \bar{w}_{k,i,j}^{(l-1)} \cdot \bar{f}_{k,j}^{(l-1)}(x) + \bar{w}_{k,i,0}^{(l-1)} \right)
\end{equation*}
for $l=2, \dots, L$
and
\begin{equation*}
\bar{f}_{k,i}^{(1)}(x) = \sigma \left(\sum_{j=1}^d \bar{w}_{k,i,j}^{(0)} \cdot x^{(j)} + \bar{w}_{k,i,0}^{(0)} \right).
\end{equation*}

In the \textit{first step} of the proof we show
\[
f_{\bar{\bw}}(x) \geq 1 - \frac{1}{n} \mbox{ for all } x \in
[u^{(1)}+\delta, v^{(1)}-\delta] \times \dots \times
[u^{(d)}+\delta, v^{(d)}-\delta]. 
\]

Let $x \in [u^{(1)}+\delta, v^{(1)}-\delta] \times \dots \times [u^{(d)}+\delta, v^{(d)}-\delta] $.
Then we get for any $i \in \{1,\dots,d\}$ by (\ref{le5eq1}), (\ref{le5eq3}) and (\ref{le5eq10})
\begin{eqnarray*}
&&\sum_{j=1}^{d} \bar{w}_{1,i,j}^{(0)}\cdot x^{(j)} + \bar{w}_{1,i,0}^{(0)}\\
&&=\sum_{j=1}^{d} (\bar{w}_{1,i,j}^{(0)}-w_{1,i,j}^{(0)})\cdot x^{(j)} + (\bar{w}_{1,i,0}^{(0)}-w_{1,i,0}^{(0)})+ \sum_{j=1}^{d} w_{1,i,j}^{(0)} \cdot x^{(j)} + w_{1,i,0}^{(0)}\\
&&\geq -d (\log n) \cdot \alpha_n- \log n + \frac{4d(\log n)^2}{\delta}(u^{(i)}+ \delta) - \frac{4d(\log n)^2}{\delta}u^{(i)} \\
&&=3d (\log n)^2 - \log n\\
&&\geq \log n 
\end{eqnarray*}

and for any $i \in \{d+1,\dots,2d\}$ by (\ref{le5eq2}), (\ref{le5eq3}) and (\ref{le5eq10})

\begin{eqnarray*}
&&\sum_{j=1}^{d} \bar{w}_{1,i,j}^{(0)}\cdot x^{(j)} + \bar{w}_{1,i,0}^{(0)}\\
&&=\sum_{j=1}^{d} (\bar{w}_{1,i,j}^{(0)}-w_{1,i,j}^{(0)})\cdot x^{(j)} + (\bar{w}_{1,i,0}^{(0)}-w_{1,i,0}^{(0)})+ \sum_{j=1}^{d} w_{1,i,j}^{(0)} \cdot x^{(j)} + w_{1,i,0}^{(0)}\\
&&\geq -d (\log n) \cdot \alpha_n- \log n - \frac{4d(\log n)^2}{\delta}(v^{(i-d)}- \delta) + \frac{4d(\log n)^2}{\delta}v^{(i-d)} \\
&&=3d (\log n)^2 - \log n\\
&&\geq \log n.
\end{eqnarray*}

This implies 

\[
\bar{f}_{1,i}^{(1)}(x) \geq \sigma(\log n) = 1- \frac{1}{n+1} \geq 1- \frac{1}{n}
\]
for any $i \in \{1,\dots,2d\}$.

Using (\ref{le5eq4})-(\ref{le5eq6}) and $|\sigma(u)| \leq 1$ for any $u \in \mathbb{R}$ we get similarly as above
\begin{eqnarray*}
&&\sum_{j=1}^{r} \bar{w}_{1,1,j}^{(1)}\cdot \bar{f}_{1,j}^{(1)}(x) + \bar{w}_{1,1,0}^{(1)}\\
&&\geq -(r+1) \log n + \sum_{j=1}^{r} w_{1,1,j}^{(1)} \cdot \bar{f}_{1,j}^{(1)}(x) + w_{1,1,0}^{(1)}\\
&&= -(r+1) \log n + \sum_{j=1}^{2d} w_{1,1,j}^{(1)} \cdot \bar{f}_{1,j}^{(1)}(x) +\sum_{j=2d+1}^{r} w_{1,1,j}^{(1)} \cdot \bar{f}_{1,j}^{(1)}(x) + w_{1,1,0}^{(1)}\\
&&\geq -(r+1) \log n + 2d \cdot 8(\log n)^2 \left(1 -\frac{1}{n}\right) - 8 (\log n)^2 \left(2d-\frac{1}{2} \right) \\
&&= -(r+1) \log n + 8 (\log n)^2 \left(\frac{1}{2} - \frac{2d}{n} \right)\\
&&\geq \log n.
\end{eqnarray*}

Therefore, we obtain 
\[
\bar{f}_{1,1}^{(2)}(x) \geq 1- \frac{1}{n}.
\]

With the same argument as above and with (\ref{le5eq7})-(\ref{le5eq9}), we can recursively conclude for $l=3,\dots,L$ that we have

\begin{eqnarray*}
&&\sum_{j=1}^{r} \bar{w}_{1,1,j}^{(l-1)}\cdot \bar{f}_{1,j}^{(l-1)}(x) + \bar{w}_{1,1,0}^{(l-1)}\\
&&\geq -(r+1) \log n + \sum_{j=1}^{r} w_{1,1,j}^{(l-1)} \cdot \bar{f}_{1,j}^{(l-1)}(x) + w_{1,1,0}^{(l-1)}\\
&&\geq -(r+1) \log n + 6(\log n)^2 \left(1 -\frac{1}{n}\right) - 3 (\log n)^2 \\
&&= -(r+1) \log n + 3 (\log n)^2 - \frac{6}{n} (\log n)^2\\
&&\geq 2 \cdot (\log n)^2 - \frac{6}{n} (\log n)^2\\
&&\geq \log n.
\end{eqnarray*}

Therefore, we obtain
\[
\bar{f}_{1,1}^{(l)}(x) \geq 1- \frac{1}{n}
\]
for $l=3,\dots, L$.

This implies

\[
f_{\bar{\bw}}(x) \geq 1 - \frac{1}{n} \mbox{ if } x \in
[u^{(1)}+\delta, v^{(1)}-\delta] \times \dots \times
[u^{(d)}+\delta, v^{(d)}-\delta].
\]

In the \textit{second step} of the proof we show
\[
f_{\bar{\bw}}(x) \leq \frac{1}{n} \mbox{ if } x^{(i)} \notin
[u^{(i)}-\delta, v^{(i)}+\delta] \mbox{ for some } i \in \{1,\dots, d\}. 
\]

Let $i \in \{1,\dots,d\}$ and assume $x^{(i)} \notin [u^{(i)}-\delta, v^{(i)}+\delta]$.
Then, we can argue similarly as above. In case $x^{(i)}< u^{(i)}-\delta$ for some $i \in \{1,\dots,d\}$ we obtain  by (\ref{le5eq1}), (\ref{le5eq3}) and (\ref{le5eq10})
\begin{eqnarray*}
&&\sum_{j=1}^{d} \bar{w}_{1,i,j}^{(0)}\cdot x^{(j)} + \bar{w}_{1,i,0}^{(0)}\\
&&\leq d (\log n) \cdot \alpha_n + \log n + \frac{4d(\log n)^2}{\delta}(u^{(i)}- \delta) - \frac{4d(\log n)^2}{\delta}u^{(i)} \\
&&\leq-3d (\log n)^2 + \log n\\
&&\leq -\log n 
\end{eqnarray*}

and in case $x^{(i)} > v^{(i)}+\delta$ for some $i \in \{1,\dots,d\}$ we obtain by (\ref{le5eq2}), (\ref{le5eq3}) and (\ref{le5eq10})

\begin{eqnarray*}
&&\sum_{j=1}^{d} \bar{w}_{1,i+d,j}^{(0)}\cdot x^{(j)} + \bar{w}_{1,i+d,0}^{(0)}\\
&&\leq d (\log n) \cdot \alpha_n+ \log n - \frac{4d(\log n)^2}{\delta}(v^{(i)}+ \delta) + \frac{4d(\log n)^2}{\delta}v^{(i)} \\
&&\leq -3d (\log n)^2 + \log n\\
&&\leq -\log n.
\end{eqnarray*}

Together with the logistic squasher it holds
\[
\bar{f}_{1,i}^{(1)}(x) \leq \sigma(-\log n) = \frac{1}{n+1} \leq \frac{1}{n}
\]
for some $i \in \{1,...,2d\}$.

Similar to above we get with (\ref{le5eq4})-(\ref{le5eq6}) and (\ref{le5eq10})
\begin{eqnarray*}
&&\sum_{j=1}^{r} \bar{w}_{1,1,j}^{(1)}\cdot \bar{f}_{1,j}^{(1)}(x) + \bar{w}_{1,1,0}^{(1)}\\
&&\leq (r+1) \log n + \sum_ {j=1}^{r} w_{1,1,j}^{(1)} \cdot \bar{f}_{1,j}^{(1)}(x) + w_{1,1,0}^{(1)}\\
&& \leq (r+1) \log n + (2d-1) \cdot 8(\log n)^2 + 8(\log n)^2 \cdot \frac{1}{n} - 8(\log n)^2 \left(2d - \frac{1}{2} \right)\\
&&=(r+1) \log n + 8(\log n)^2 \left(\frac{1}{n} - \frac{1}{2} \right)\\
&&\leq - \log n.
\end{eqnarray*}

From this we obtain
\[
\bar{f}_{1,1}^{(2)}(x) \leq \frac{1}{n}.
\]

Using (\ref{le5eq7})-(\ref{le5eq10}) we can argue similarly as above and conclude recursively for $l=3,\dots,L$
\begin{eqnarray*}
&&\sum_{j=1}^{r} \bar{w}_{1,1,j}^{(l-1)}\cdot \bar{f}_{1,j}^{(l-1)}(x) + \bar{w}_{1,1,0}^{(l-1)}\\
&&\leq (r+1) \log n + \sum_{j=1}^{r}w_{1,1,j}^{(l-1)} \cdot \bar{f}_{1,j}^{(l-1)}(x)+w_{1,1,0}^{(l-1)}\\
&&\leq (r+1) \log n + 6 \cdot (\log n)^2 \cdot \frac{1}{n} - 3 (\log n)^2\\
&&\leq (\log n)^2 + 6 \cdot (\log n)^2 \cdot \frac{1}{n} - 3 (\log n)^2\\
&&= -2 \cdot (\log n)^2 + 6 \cdot (\log n)^2 \cdot \frac{1}{n}\\
&&\leq - \log n.
\end{eqnarray*}

So we get
$
\bar{f}_{1,1}^{(l)}(x) \leq \frac{1}{n}
$
for $l=3,\dots, L$.
Therefore,
\[
f_{\bar{\bw}}(x) \leq \frac{1}{n} \mbox{ if } x^{(i)} \notin
[u^{(i)}-\delta, v^{(i)}+\delta] \mbox{ for some } i \in \{1,\dots, d\}
\]
holds.

This yields the assertion.
\quad  \hfill $\Box$

\begin{lemma}
  \label{le6}
  Let $0 < \delta \leq 1 $, $1 \leq \alpha_n \leq \log n$ and let $\sigma$ be the logistic squasher,
  let $m:\Rd \rightarrow \R$ be Lipschitz continuous with Lipschitz
  constant $C_{Lip}$, let $L, r, n \in \N$ with $L \geq 2$, $r \geq 2d$, $n \geq 8d$, $n \geq\exp(r+1)$
  and let $K \in \N$ with $K^{d} \leq K_n$. Furthermore define $f_{\bar{\bw}}$ by 
\begin{equation*}
f_{\bar{\bw}}(x) = \sum_{j=1}^{K_n} \bar{w}_{1,1,j}^{(L)} \cdot \bar{f}_{j,1}^{(L)}(x) 
\end{equation*}  
for some $\bar{w}_{1,1,1}^{(L)}, \dots, \bar{w}_{1,1,K_n}^{(L)} \in \mathbb{R}$, where $\bar{f}_{j,1}^{(L)}$ are recursively defined by
\begin{equation*}
\bar{f}_{k,i}^{(l)}(x) = \sigma\left(\sum_{j=1}^{r} \bar{w}_{k,i,j}^{(l-1)} \cdot \bar{f}_{k,j}^{(l-1)}(x) + \bar{w}_{k,i,0}^{(l-1)} \right)
\end{equation*}
for some $\bar{w}_{k,i,0}^{(l-1)}, \dots, \bar{w}_{k,i, r}^{(l-1)} \in \mathbb{R}$
$(l=2, \dots, L)$
and
\begin{equation*}
\bar{f}_{k,i}^{(1)}(x) = \sigma \left(\sum_{j=1}^d \bar{w}_{k,i,j}^{(0)} \cdot x^{(j)} + \bar{w}_{k,i,0}^{(0)} \right).
\end{equation*}
Choose $\bw$ such that 
\begin{equation}
	\label{le6beq4}
	w_{j_k,j,j}^{(0)}=\frac{4d \cdot (\log n)^2}{\delta}
	\quad \mbox{and} \quad
	w_{j_k,j,0}^{(0)}=\frac{-4d \cdot (\log n)^2 \cdot u_k^{(j)}}{\delta}
	\quad \mbox{for } j \in \{1, \dots, d\},
\end{equation}
\begin{equation}
	\label{le6beq5}
	w_{j_k,j+d,j}^{(0)}=\frac{-4d \cdot (\log n)^2}{\delta}
	\quad \mbox{and} \quad
	w_{j_k,j+d,0}^{(0)}=\frac{4d \cdot (\log n)^2 \cdot v_k^{(j)}}{\delta}
	\quad \mbox{for }j \in \{1, \dots, d\},
\end{equation}
\begin{equation}
	\label{le6beq6}
	w_{j_k,s,t}^{(0)}=0
	\quad \mbox{if }
	s \leq 2d,
	s \neq t, s \neq t+d \mbox{ and } t>0,
\end{equation}
\begin{equation}
	\label{le6beq7}
	w_{j_k,1,t}^{(1)}= 8 \cdot (\log n)^2\quad \mbox{for } t \in \{1, \dots, 2d\},
\end{equation}
\begin{equation}
	\label{le6beq8}
	w_{j_k,1,0}^{(1)} = - 8 \cdot (\log n)^2 \left(2d-\frac{1}{2}\right)
\end{equation}
\begin{equation}
	\label{le6beq9}
	w_{j_k,1,t}^{(1)}=0
	\quad \mbox{for }
	t > 2d,
\end{equation}
\begin{equation}
	\label{le6beq10}
	w_{j_k,1,1}^{(l)}= 6 \cdot (\log n)^2 \quad \mbox{for } l \in \{2, \dots, L\},
\end{equation}
\begin{equation}
	\label{le6beq11}
	w_{j_k,1,0}^{(1)} = -3 (\log n)^2 \quad \mbox{for }l \in \{2, \dots, L\}
\end{equation}
and
\begin{equation}
	\label{le6beq12}
	w_{j_k,1,t}^{(l)}=0
	\quad \mbox{for }
	t > 1 \mbox{ and } l \in \{2, \dots, L\}
\end{equation} 
for all $k \in \{1, \dots, {K^{d}}\}$.

Let $a_{1},\dots, a_{d},b_{1}, \dots, b_{d}  \in [-\alpha_n, \alpha_n]^d$
  with $b_{i}-a_{i}=\Delta$ for all $i \in \{1,\dots, d\}$ and $\Delta \in \mathbb{R}_{+}$. Then there exist
  \begin{equation}
    \label{le6eq1}
    \alpha_1, \dots, \alpha_{K^{d}} \in [- \|m\|_\infty,
  \|m\|_\infty]
  \end{equation}
  and
  $u_1,v_1, \dots, u_{K^{d}},v_{K^{d}} \in [a_{1},b_{1})\times \dots \times [a_{d},b_{d})$
  such that for all
  pairwise distinct $j_1, \dots, j_{K^{d}} \in \{1, \dots, K_n\}$ the inequality
\begin{equation}\label{le6eq2}
  | f_{\bar{\bw}}(x)-m(x)|  \leq
  c_{19} \cdot \left(
  C_{Lip} \cdot
  \frac{\Delta}{K} + K^{d} \cdot \frac{1}{n} 
  \right)
 \end{equation}
  holds for all $x \in [a_{1},b_{1}]\times \dots \times [a_{d},b_{d}]$ which are not contained in
\begin{equation}\label{set}
  \bigcup_{j \in \{0,1,...,K\}} \bigcup_{i \in \{1,...,d\}} \left\{x \in \mathbb{R}^d:\left|x^{(i)}- \left(
a_i + j \cdot \frac{b_i-a_i}{K}
  \right)\right| < \delta \right\} 
\end{equation}  
and for all weight vectors $\bar{\bw}$ which satisfy
  \begin{equation}
    \label{le6beq2}
    \bar{w}_{1,1,j_k}^{(L)}=\alpha_k \quad (k \in \{1, \dots, K^{d}\}), \quad
    \bar{w}_{1,1,k}^{(L)}=0 \quad (k \notin \{j_1, \dots, j_{K^{d}}\})
    \end{equation}
    and
  \begin{equation}
    \label{le6beq3}
    |w_{j_s,k,i}^{(l)}-\bar{w}_{j_s,k,i}^{(l)}| \leq \log n
    \quad \mbox{for all } l \in \{0, \dots, L-1\}, s \in \{1, \dots, K^{d}\}.
  \end{equation}

For $\delta \leq \frac{\Delta}{K}$ and $x \in \mathbb{R}^d$ we get additionally
\begin{equation}\label{le6eq13}
| f_{\bar{\bw}}(x)|  \leq \|m\|_\infty \cdot \left( 3^d + \frac{K^{d}}{n} \right).
\end{equation}
  \end{lemma}

\noindent
{\bf Proof.}
We partition $[a_{1},b_{1})\times \dots \times [a_{d},b_{d})$ into $K^{d}$ equivolume cubes of side length $\frac{\Delta}{K}$.
For comprehensibility, we number these cubes $C_i$ by $i \in \{1,\dots,K^{d}\}$, such that $C_i$ corresponds to the cube
\begin{equation*}
[u_{i}^{(1)},v_{i}^{(1)}) \times \dots \times [u_{i}^{(d)},v_{i}^{(d)}).
\end{equation*}

Since $m$ is Lipschitz continuous with Lipschitz constant $C_{Lip}$, $m$ can be approximated by an approximand $S$ that is piecewise constant on each cube. $S$ can be expressed in the form
\begin{equation*}
S(x) = \sum_{i \in \{1,\dots,K^d\}} \alpha_i \cdot 1_{C_i}(x)
\end{equation*}
where $\alpha_i = m(z_i)$ for $i\in\{1,\dots,K^{d}\}$ and $z_i$ is the center of the cube $C_{i}$.

$S(x)$ has value $m(z_i$) for $x \in C_i$. Therefore it follows from the Lipschitz continuity of $m$
\[
|S(x)-m(x)| \leq C_{Lip} \cdot \|z_i-x\|_2
\]
if $x \in C_i$ for some $i\in\{1,\dots,K^{d}\}$.
Since every cube has a side length of $\frac{\Delta}{K}$ we obtain

\[
|S(x)-m(x)| \leq C_{Lip} \cdot \sqrt{d} \cdot \frac{\Delta}{K}
\mbox{ for }
x \in [a_{1},b_{1})\times \dots \times [a_{d},b_{d}).
\]
Now, because of (\ref{le6beq2}) and the definitions of $S$ and $f_{\bar{\bw}}(x)$, we have
\begin{align*}
S(x)-f_{\bar{\bw}}(x)= \sum_{k=1}^{K^d} \alpha_k \left(1_{C_k}(x) - \bar{f}_{j_k,1}^{(L)}(x) \right).
\end{align*}
Application of Lemma \ref{le5} yields
\begin{equation*}
|S(x)-f_{\bar{\bw}}(x)| \leq c_{21} \cdot K^{d} \cdot \frac{1}{n}
\end{equation*}
for all $x \in [a_{1},b_{1})\times \dots \times [a_{d},b_{d})$ which are not contained in (\ref{set}).

From this we obtain
\[
\big| f_{\bar{\bw}}(x)-m(x)\big|= \big| f_{\bar{\bw}}(x)-S(x)\big| +\big|S(x)-m(x)\big|
\leq c_{22} \left( K^{d} \cdot \frac{1}{n} + C_{Lip} \cdot \frac{\Delta}{K}\right) 
\]
for all $x \in [a_{1},b_{1})\times \dots \times [a_{d},b_{d})$ which are not contained in (\ref{set}).
This implies (\ref{le6eq2}).

To show inequality (\ref{le6eq13}), we assume that $x \in \mathbb{R}^d$.
Then for $\delta \leq \frac{\Delta}{K}$ and every fixed $x$ we get
\begin{eqnarray*}
|f_{\bar{\bw}}(x)| &=& \sum_{k=1}^{K^{d}} \bar{w}_{1,1,j_k}^{(L)} \cdot \bar{f}_{j_k,1}^{(L)}(x)\\
  &=& \sum_{k \in \{1, \dots, K^d \}\ :\  x \in
[u_{k}^{(1)}-\delta,v_{k}^{(1)}+\delta] \times \dots \times [u_{k}^{(d)}-\delta,v_{k}^{(d)}+\delta]
} \bar{w}_{1,1,j_k}^{(L)} \cdot \bar{f}_{j_k,1}^{(L)}(x) \\
  && 
  + \sum_{k \in \{1, \dots, K^d\}\ :\  x \not\in
    [u_{k}^{(1)}-\delta,v_{k}^{(1)}+\delta] \times \dots \times [u_{k}^{(d)}-\delta,v_{k}^{(d)}+\delta]
  } \bar{w}_{1,1,j_k}^{(L)}\cdot \bar{f}_{j_k,1}^{(L)}(x)\\
\end{eqnarray*}
There are at most $3^d$ many cubes $[u_{i}^{(1)}-\delta,v_{i}^{(1)}+\delta] \times \dots \times [u_{i}^{(d)}-\delta,v_{i}^{(d)}+\delta]$ which contain $x$. Together with the definition of $\bar{w}_{1,1,j_k}$ for $k\in \{1, \dots, K^{d}\}$ we get
\[
\sum_{k\in \{1, \dots, K^d\}\ :\  x \in
[u_{k}^{(1)}-\delta,v_{k}^{(1)}+\delta] \times \dots \times [u_{k}^{(d)}-\delta,v_{k}^{(d)}+\delta] } \bar{w}_{1,1,j_k}^{(L)} \cdot \bar{f}_{j_k,1}^{(L)}(x)
\leq 3^d \cdot \|m\|_\infty.
\]

By Lemma \ref{le5}
we know that $\bar{f}_{j_k,1}^{(L)}(x) \leq \frac{1}{n}$ for $k \in \{1,\dots,K^{d}\}$ if $x \not\in
[u_{k}^{(1)}-\delta,v_{k}^{(1)}+\delta] \times \dots \times [u_{k}^{(d)}-\delta,v_{k}^{(d)}+\delta]$. Thus we get together with the definition of  $w_{1,1,j_k}$ for $k\in \{1, \dots, K^{d}\}$
\[\sum_{k\in \{1, \dots, K^d\}\ :\  x \not\in
[u_{k}^{(1)}-\delta,v_{k}^{(1)}+\delta] \times \dots \times [u_{k}^{(d)}-\delta,v_{k}^{(d)}+\delta]
  } w_{1,1,j_k}^{(L)}\cdot f_{j_k,1}^{(L)}(x) \leq K^{d} \cdot \frac{1}{n} \cdot \|m\|_\infty.\]

This results in
\[
|f_{\bar{\bw}}(x)| \leq  3^d \cdot \|m\|_\infty + K^{d} \cdot \frac{1}{n} \cdot \|m\|_\infty
\]
for $x \in \mathbb{R}^d$.

\quad  \hfill $\Box$

\begin{lemma}\label{le7}

Let $\sigma$ be the logistic squasher, let $1 \leq \alpha_n \leq \log n$,
let $m:\Rd \rightarrow \R$ be Lipschitz continuous as well as bounded,
let $L,r,n \in \N$ with $L \geq 2$, $r \geq 2d$, $n \geq 8d$ and $n \geq\exp(r+1)$
and let $K \in \N$ with
$ 2 \leq K \leq \alpha_n -1$ and
$(K^2+1)^{3d} \leq K_n$. Choose $\bw$ such that 
\begin{equation}
	\label{le7eq4}
	w_{j_k,j,j}^{(0)}=4d \cdot K^2 \cdot (\log n)^2
	\quad \mbox{and} \quad w_{j_k,j,0}^{(0)}=-4d \cdot K^2 \cdot (\log n)^2 \cdot u_k^{(j)}
	\quad \mbox{for } j \in \{1, \dots, d\},
\end{equation}
\begin{equation}
	\label{le7eq5}
	w_{j_k,j+d,j}^{(0)}=-4d \cdot K^2 \cdot (\log n)^2
	\quad \mbox{and} \quad
	w_{j_k,j+d,0}^{(0)}=4d \cdot K^2 \cdot (\log n)^2 \cdot v_k^{(j)}
	\quad\mbox{for } j \in \{1, \dots, d\},
\end{equation}
\begin{equation}
	\label{le7eq6}
	w_{j_k,s,t}^{(0)}=0
	\quad \mbox{if }
	s \leq 2d,
	s \neq t, s \neq t+d \mbox{ and } t>0,
\end{equation}
\begin{equation}
	\label{le7eq7}
	w_{j_k,1,t}^{(1)}= 8 \cdot (\log n)^2\quad \mbox{for } t \in \{1, \dots, 2d\}),
\end{equation}
\begin{equation}
	\label{le7eq8}
	w_{j_k,1,0}^{(1)} = - 8 \cdot (\log n)^2 \left(2d-\frac{1}{n}\right)
\end{equation}
\begin{equation}
	\label{le6beq9}
	w_{j_k,1,t}^{(1)}=0
	\quad \mbox{for }
	t > 2d,
\end{equation}
\begin{equation}
	\label{le7eq10}
	w_{j_k,1,1}^{(l)}= 6 \cdot (\log n)^2 \quad \mbox{for } l \in \{2, \dots, L\},
\end{equation}
\begin{equation}
	\label{le7eq11}
	w_{j_k,1,0}^{(l)} = -3 \cdot (\log n)^2 \quad \mbox{for } l \in \{2, \dots, L\}
\end{equation}
and
\begin{equation}
	\label{le7eq12}
	w_{j_k,1,t}^{(l)}=0
	\quad \mbox{for }
	t > 1 \mbox{ and } l \in \{2, \dots, L\}
\end{equation} 
 for all $k \in \{1, \dots, (K^2+1)^{3d}\}$.

Then there exists
\begin{equation}
  \label{le7eq1}
  \alpha_1, \dots, \alpha_{(K^2+1)^{3d}} \in
        \left[- \frac{\|m\|_\infty}{(K^2+1)^{2d}},
\frac{\|m\|_\infty}{(K^2+1)^{2d}} \right]
\end{equation}
and
$u_1,v_1, \dots, u_{(K^2+1)^{3d}},v_{(K^2+1)^{3d}} \in [-K-\frac{2}{K},K]^d$
such that for all
pairwise distinct\\ $j_1, \dots, j_{(K^2+1)^{3d}} \in \{1, \dots, K_n\}$  
\begin{eqnarray}\label{le7eq1}
&&\int \left|f_{\bar{\bw}}(x)-m(x)\right|^2 \PROB_X(dx)\nonumber\\
  &&\leq
  c_{23} \cdot \left( \frac{1}{K} + \frac{K^{12d}}{n^2}
  +\left(\frac{K^{6d}}{n} +1\right)^2 \cdot \PROB_X(\mathbb{R}^d\setminus [-K,K]^d)\right)
\end{eqnarray}
holds
for all weight vectors $\bar{\bw}$ which satisfy 
 \begin{equation}
	\label{le7beq2}
	\bar{w}_{1,1,j_k}^{(L)}=\alpha_k \quad (k \in \{1, \dots, (K^2+1)^{3d}\}), \quad
	\bar{w}_{1,1,k}^{(L)}=0 \quad (k \notin \{j_1, \dots, j_{(K^2+1)^{3d}}\})
\end{equation}
and
\begin{equation}
	\label{le7eq3}
	|w_{j_s,k,i}^{(l)}-\bar{w}_{j_s,k,i}^{(l)}| \leq \log n
	\quad \mbox{for all } l \in \{0, \dots, L-1\}, s \in \{1, \dots, (K^2+1)^{3d}\}.
\end{equation}

\end{lemma}

\noindent
{\bf Proof.}
We subdivide $[-K-\frac{2}{K},K]^d$ in $(K^2+1)^{d}$ cubes of side length $\frac{2}{K}$.
We number these cubes $C_{i}$ by $i \in \{1,\dots,(K^2+1)^{d}\}$, such that $C_{i}$ corresponds to the cube

\[
[u_{i}^{(1)},v_{i}^{(1)}) \times \dots \times [u_{i}^{(d)},v_{i}^{(d)}).
    \]
    Let $C_{Lip}$ be the Lipschitz constant of $m$.
Then by Lemma \ref{le6} applied to $m/(K^2+1)^{2d}$ and $\delta=\frac{1}{K^2}$
we know that
\[
\left| f_{\bar{\bw}}(x)- \frac{1}{(K^2+1)^{2d}} \cdot m(x)
\right|  \leq
  c_{19} \cdot \left(
   \frac{C_{Lip}}{(K^2+1)^{2d}} \cdot
  \frac{2}{K} + (K^2+1)^{d} \cdot \frac{1}{n} 
  \right)
\]
holds for all $x \in [-K-\frac{2}{K},K]^d$ which are not contained in
\begin{equation}\label{le7eq2}
  A :=
    \bigcup_{i \in \{0,1,...,K+1\}} \bigcup_{j \in \{1,...,d\}} \left\{x \in \mathbb{R}^d:\left|x^{(j)}- \left(
-K - \frac{2}{K} + i \cdot \frac{2}{K}
  \right)\right| < \delta \right\}. 
\end{equation}
Next we repeat the whole construction $(K^2+1)^{2d}$ many times, which results
in an approximation   $f_{\bar{\bw}}$ of
\[
(K^2+1)^{2d} \cdot  \frac{1}{(K^2+1)^{2d}} \cdot m(x)
\]
which satisfies
\begin{equation}
  \label{ple7eq*}
\left| f_{\bar{\bw}}(x)-  m(x)
\right|  \leq
  c_{24} \cdot \left(
  \frac{1}{K} + (K^2+1)^{3d} \cdot \frac{1}{n} 
  \right)
  \end{equation}
  outside of $A$.

Now we want to move the grid so that $[-K,K]^d$ is always covered.
We slightly shift the whole grid of cubes along the $j$-th component by modifying all $u_{i}^{(j)},v_{i}^{(j)}$ by the same additional summand. This summand is
chosen from the set
\[
\left\{
k \cdot
\frac{2}{K^2} \quad : \quad k=0,1, \dots, K-1
\right\}
\]
for fixed $j \in \{1,\dots,d\}$. In this way we can construct $K$
different versions of $f_{\bar{\bw}}$ that still satisfy (\ref{ple7eq*}) for all $x \in [-K,K]^d$ up to corresponding versions of $A$.

Since we shift the grid of cubes we obtain for fixed $j \in \{1,\dots, d\}$ $K$ disjoint versions of $\bigcup_{i \in \{0,1,...,K+1\}} \left\{x \in \mathbb{R}^d:\left|x^{(j)}- \left(
-K - \frac{2}{K} + i \cdot \frac{2}{K}
\right)\right| < \delta \right\}$.
Because the sum of $\PROB_X$-measures of these $K$ disjoint sets is less than or equal to one, at least one of them must have measure less than or equal to $\frac{1}{K}$. Consequently we can shift the $u_i$ and $v_i$ such that
\begin{equation*}
\PROB_X \left(A \right)
= \sum_{j \in \{1,\dots,d\}} \frac{1}{K} = \frac{d}{K}
\end{equation*}
holds.

Now we have found a shifted version of the grid such that the set (\ref{le7eq2}) has a measure less than or equal to $\frac{d}{K}$.
By  (\ref{ple7eq*}) we know that   $|f_{\bar{\bw}}(x)-m(x)| \leq c_{25} \cdot \left(\ \frac{1}{K} + \frac{K^{6d}}{n} \right)$ holds for $x \in [-K,K]^d\setminus A$. From this together with the second assertion from Lemma \ref{le6} we obtain

\begin{eqnarray*}
&&\int \left|f_{\bar{w}}(x)-m(x)\right|^2 \PROB_X(dx)\\
&&=\int_{[-K,K]^d\setminus A}\left|f_{\bar{w}}(x)-m(x)\right|^2 \PROB_X(dx) +
\int_{A}\left|f_{\bar{w}}(x)-m(x)\right|^2 \PROB_X(dx)\\
&& \qquad + \int_{\mathbb{R}^d\setminus {[-K,K]^d}} \left|f_{\bar{w}}(x)-m(x)\right|^2 \PROB_X(dx)\\
&&\leq c_{25}^2 \left(\frac{1}{K} + \frac{K^{6d}}{n}\right)^2  + c_{26}
\left(3^d + \frac{K^{6d}}{n}\right)^2  \cdot \frac{d}{K}
\\
&& \quad
+ c_{27} \left(3^d +\frac{ K^{6d}}{n} \right)^2
\cdot \PROB_X(\mathbb{R}^d\setminus [-K,K]^d),
\end{eqnarray*}
which implies the assertion.
\quad  \hfill $\Box$

In order to be able to formulate our next auxiliary result we
need the following notation:
Let $(x_1,y_1), \dots, (x_n,y_n) \in \Rd \times \R$, let $K \in \N$,
let $B_1,\dots,B_K:\Rd \rightarrow \R$ and let $c_2>0$. In the next lemma
we consider the problem to minimize
\begin{equation}
  \label{se5eq1}
  F(\ba) =
  \frac{1}{n} \sum_{i=1}^n
  |\sum_{k=1}^K a_k \cdot B_k(x_i)-y_i|^2
  +
    c_2 \cdot  \sum_{k=1}^{K_n} a_k^2 ,
  \end{equation}
where $\ba=(a_1,\dots,a_K)^T$,
by gradient descent. To do this, we choose $\ba^{(0)} \in \R^K$
and set
\begin{equation}
  \label{se5eq2}
  \ba^{(t+1)} = \ba^{(t)}
  - \lambda_n \cdot (\nabla_\ba F)(\ba^{(t)})
    \end{equation}
for some properly chosen $\lambda_n>0.$

        \begin{lemma}
          \label{le8}
          Let $F$ be defined by (\ref{se5eq1}) and choose $\ba_{opt}$
          such that
          \[
F(\ba_{opt})=\min_{\ba \in \R^{K}} F(\ba).
          \]
          Then for any
          $\ba \in \R^{K}$
          we have
                    \[
                            \|(\nabla_\ba F)(\ba)\|^2
                            \geq 4 \cdot c_2  \cdot (F(\ba)-F(\ba_{opt})).
                            \]
          \end{lemma}

        \noindent
            {\bf Proof.}
            The proof is a modification of the proof of Lemma 3 in Braun, Kohler and Walk (2019).

            Set
            \[
            \bE= c_2 \cdot \left(
            \begin{array}{ccccc}
              1 & 0 & 0 & \dots & 0 \\
              0 & 1 & 0 & \dots & 0 \\
              \vdots & \vdots & \vdots & \vdots & \vdots \\
              0 & 0 & 0 & \dots & 1 \\
              \end{array}
            \right),
            \]
            \[
\bB = \left(
B_j(x_i)
\right)_{1 \leq i \leq n, 1 \leq j \leq K}
\quad \mbox{and} \quad
\bA=\frac{1}{n} \cdot \bB^T \cdot \bB + c_2 \cdot \bE.
            \]
            Then $\bA$ is positive definite
            and hence regular, from which we can conclude
            \begin{eqnarray*}
              F(\ba)         & = &
              \frac{1}{n} \cdot \left( \bB \cdot \ba - \by \right)^T \cdot \left( \bB \cdot \ba - \by \right) +
c_2 \cdot \ba^T \cdot \bE \cdot \ba
              \\
              &=&
              \ba^T \bA \ba - 2 \by^T \frac{1}{n} \bB \ba + \frac{1}{n} \by^T \by
              \\
              &=&
              (\ba - \bA^{-1} \frac{1}{n} \bB^T \by)^T \bA (\ba - \bA^{-1} \frac{1}{n} \bB^T \by)
              + F(\ba_{opt}),
              \end{eqnarray*}
            where
            \[
            F(\ba_{opt})
            =
            \frac{1}{n} \by^T \by
            -
            \by^T \cdot \frac{1}{n} \cdot \bB \cdot \bA^{-1} \cdot \frac{1}{n} \cdot
            \bB^T \by.
            \]
            Using
            \[
            \bb^T \bA \bb \geq
            c_2 \cdot \bb^T \bE \bb
            =
            c_2 \cdot \bb^T \bb
            \]
   and $\bA^T=\bA$         we conclude
   \begin{eqnarray*}
     &&
     F(\ba)-F(\ba_{opt})
     \\
     &&
     =
     ((\bA^{1/2})^T  (\ba - \bA^{-1} \frac{1}{n} \bB^T \by))^T \bA^{1/2} (\ba - \bA^{-1} \frac{1}{n} \bB^T \by)
     \\
     &&
     \leq
     \frac{1}{c_2}
     \cdot
     ((\bA^{1/2})^T  (\ba - \bA^{-1} \frac{1}{n} \bB^T \by))^T \bA \bA^{1/2} (\ba - \bA^{-1} \frac{1}{n} \bB^T \by)
     \\
     &&
     =
    \frac{1}{c_2}
     \cdot
          ((\bA)^T  (\ba - \bA^{-1} \frac{1}{n} \bB^T \by))^T \bA  (\ba - \bA^{-1} \frac{1}{n} \bB^T \by)
  \\
     &&
     =
    \frac{1}{c_2}
     \cdot
          (\bA  \ba - \frac{1}{n} \bB^T \by)^T (\bA \ba -  \frac{1}{n} \bB^T \by)\\
     &&
     =
    \frac{1}{4 \cdot c_2}
     \cdot
          (2 \bA  \ba - \frac{2}{n} \bB^T \by)^T (2 \bA \ba -  \frac{2}{n} \bB^T \by)\\
     &&
     =
    \frac{1}{4 \cdot c_2}
     \cdot
   \left\| (\nabla_\ba F)(\ba) \right\|^2 ,
   \end{eqnarray*}
   where the last equality follows from
   \begin{eqnarray*}
     &&
     (\nabla_\ba F)(\ba)
     =
     \nabla_\ba \left(
              \ba^T \bA \ba - 2 \by^T \frac{1}{n} \bB \ba + \frac{1}{n} \by^T \by
              \right)
              =
2 \bA  \ba - \frac{2}{n} \bB^T \by.
     \end{eqnarray*}
            \quad
            \hfill $\Box$

            \subsection{Proof of Theorem \ref{th1}}
            Let $\epsilon >0$ and $K \in \N$ be arbitrary.
            W.l.o.g. we assume $K_n \geq (K^2+1)^{3d}$.
            Choose a Lipschitz continuous and
            bounded function $\bar{m}:\Rd \rightarrow \R$
            such that
            \begin{equation}
              \label{pth1eq0}
\int | \bar{m}(x)-m(x)|^2 \PROB_X (dx) \leq \epsilon.
            \end{equation}
            Let $A_n$ be the event that firstly the weight vector $\bw^{(0)}$
            satisfies
            \[
            | (\bw^{(0)})_{j_s,k,i}^{(l)}-\bw_{j_s,k,i}^{(l)}| \leq \log n
            \quad \mbox{for all } l \in \{0, \dots, L-1\},
            s \in \{1, \dots, (K^2+1)^{3d} \}
            \]
            for some weight vector $\bw$ which satisfies
            the conditions
(\ref{le7eq4})--(\ref{le7eq12})
            of Lemma \ref{le7} for $\bar{m}$
            and some $j_1, \dots, j_{(K^2+1)^{3d}} \in \{1, \dots, K_n\}$,
            and that secondly
\[
\frac{1}{n} \sum_{i=1}^n Y_i^2 \leq \beta_n^3
\]
holds.
Define the weight vectors
$(\bw^*)^{(t)}$ by
\[
((\bw^*)^{(t)})_{k,i,j}^{(l)} = (\bw^{(t)})_{k,i,j}^{(l)} \quad
\mbox{for all } l=0,\dots, L-1
\]
and
\[
((\bw^*)^{(t)})_{1,1,j_k}^{(L)} = \alpha_k \quad \mbox{for all } k=1,\dots, (K^2+1)^{3d}
\]
and
\[
((\bw^*)^{(t)})_{1,1,k}^{(L)} = 0 \quad \mbox{for all } k \notin \{j_1,\dots,
j_{(K^2+1)^{3d}}\}
\]
where $\alpha_k$ is chosen as in Lemma \ref{le7} in case that $A_n$ holds and where
$\alpha_1=...=\alpha_{(K^2+1)^{3d}}=0$ in case that $A_n$ does not hold.

In the {\it first step of the proof} we 
decompose the  $L_2$ error of $m_n$ in a sum of several terms.
We have
\begin{eqnarray*}
  &&
    \int | m_n(x)-m(x)|^2 \PROB_X (dx)
    \\
    &&
    =
    (\EXP\{ |m_n(X)-Y|^2 | \D_n\} - \EXP\{ |m(X)-Y|^2 \}) \cdot 1_{A_n}
    +
    \int | m_n(x)-m(x)|^2 \PROB_X (dx) \cdot 1_{A_n^c}
    \\
    &&
    =
    (\EXP\{ |m_n(X)-Y|^2 | \D_n\}
        \\
    &&
    \hspace*{3cm}
    -
    \EXP\{ |m_n(X)-Y|^2 \cdot 1_{[-\alpha_n,\alpha_n]^d}(X) | \D_n\}) \cdot 1_{A_n}
    \\
    &&
   \quad +
(    \EXP\{ |m_n(X)-Y|^2 \cdot 1_{[-\alpha_n,\alpha_n]^d}(X) | \D_n\}
        \\
    &&
    \hspace*{3cm}
    -
    (1+\epsilon) \cdot
    \EXP\{ |m_n(X)-T_{\beta_n} Y|^2 \cdot 1_{[-\alpha_n,\alpha_n]^d}(X)| \D_n\})
    \cdot 1_{A_n}
    \\
    &&
    \quad
    +
    ((1+\epsilon) \cdot
    \EXP\{ |m_n(X)-T_{\beta_n} Y|^2 \cdot 1_{[-\alpha_n,\alpha_n]^d}(X)| \D_n\}
    \\
    &&
    \hspace*{3cm}
    -
    (1+\epsilon) \cdot \frac{1}{n} \sum_{i=1}^n |m_n(X_i)- T_{\beta_n} Y_i|^2
    \cdot 1_{[-\alpha_n,\alpha_n]^d}(X_i) )\cdot 1_{A_n}
    \\
    &&
     \quad
    +
    ((1+\epsilon) \cdot \frac{1}{n} \sum_{i=1}^n |m_n(X_i)- T_{\beta_n} Y_i|^2
    \cdot 1_{[-\alpha_n,\alpha_n]^d}(X_i)
    \\
    &&
    \hspace*{3cm}
    -
    (1+\epsilon) \cdot \frac{1}{n} \sum_{i=1}^n |f_{\bw^{(t_n)}}(X_i) - T_{\beta_n} Y_i|^2 \cdot 1_{[-\alpha_n,\alpha_n]^d}(X_i)) \cdot 1_{A_n}
    \\
    &&
    \quad
    +
    (1+\epsilon) \cdot (\frac{1}{n} \sum_{i=1}^n |f_{\bw^{(t_n)}}(X_i) - T_{\beta_n} Y_i|^2 \cdot 1_{[-\alpha_n,\alpha_n]^d}(X_i)
    \\
    &&
    \hspace*{3cm}
    -
    (1+\epsilon) \cdot \frac{1}{n} \sum_{i=1}^n |f_{\bw^{(t_n)}}(X_i) - Y_i|^2 \cdot 1_{[-\alpha_n,\alpha_n]^d}(X_i))\cdot 1_{A_n}
   \\
    &&
   \quad
   +
   ((1+\epsilon)^2 \cdot \frac{1}{n} \sum_{i=1}^n |f_{\bw^{(t_n)}}(X_i) - Y_i|^2 \cdot 1_{[-\alpha_n,\alpha_n]^d}(X_i)
   -
   \EXP\{ |m(X)-Y|^2 \}) \cdot 1_{A_n}\\
   &&
   \quad
   +
    \int | m_n(x)-m(x)|^2 \PROB_X (dx) \cdot 1_{A_n^c}   
   \\
   &&
   =
   \sum_{j=1}^7 T_{j,n}.
  \end{eqnarray*}

In the {\it second step of the proof} we show
\[
\limsup_{n \rightarrow \infty} \EXP T_{j,n} \leq 0 \quad
\mbox{for } j \in \{1,2,5\}.
\]
Because of $\alpha_n \rightarrow \infty$ $(n \rightarrow \infty)$
and
$\EXP Y^2 < \infty$ it holds
\[
\EXP T_{1,n} = \EXP\{ |Y|^2 \cdot 1_{\R^d \setminus [-\alpha_n,\alpha_n]^d}(X) \}
\rightarrow 0 \quad (n \rightarrow \infty).
\]
Using $(a+b)^2 \leq (1+\epsilon) \cdot a^2 + (1 + \frac{1}{\epsilon}) \cdot b^2$
$(a,b \in \R)$ we get
\[
\EXP T_{2,n} \leq \left(1 + \frac{1}{\epsilon} \right) \cdot \EXP\{|T_{\beta_n}Y-Y|^2\}
\]
and
\[
\EXP T_{5,n} \leq (1+\epsilon) \cdot \left(1 + \frac{1}{\epsilon} \right) \cdot \EXP\{|T_{\beta_n}Y-Y|^2\}.
\]
Because of $\beta_n \rightarrow \infty$ $(n \rightarrow \infty)$ and
$\EXP Y^2 < \infty$ this implies the assertion of the second step.

In the {\it third step of the proof} we show
\[
\limsup_{n \rightarrow \infty} \EXP T_{4,n} \leq 0. 
\]
If $|y| \leq \beta_n$ then it holds for
any $z \in \R$
\[
| T_{\beta_n} z - y| \leq |z-y|.
\]
This implies
\begin{eqnarray*}
  &&
  \frac{1}{n} \sum_{i=1}^n |m_n(X_i)- T_{\beta_n} Y_i|^2     \cdot 1_{[-\alpha_n,\alpha_n]^d}(X_i)
  \\
&&=
\frac{1}{n} \sum_{i=1}^n |T_{\beta_n}f_{\bw^{(t_n)}}(X_i) - T_{\beta_n} Y_i|^2     \cdot 1_{[-\alpha_n,\alpha_n]^d}(X_i)
\\
&&\leq
\frac{1}{n} \sum_{i=1}^n |f_{\bw^{(t_n)}}(X_i) - T_{\beta_n} Y_i|^2
\cdot 1_{[-\alpha_n,\alpha_n]^d}(X_i),
\end{eqnarray*}
hence $T_{4,n} \leq 0$ holds, which implies the
assertion of the third step.

In the {\it fourth step of the proof} we show that the
assumptions of Lemma \ref{le1} are satisfied if $A_n$ holds.
If $A_n$ holds, then we have
\[
F_n(\bw^{(0)})
  =
  \frac{1}{n} \sum_{i=1}^n Y_i^2 \leq \beta_n^3.
\]
Hence if  the assumptions of the two conditions, which we have
to show in Lemma \ref{le1}, hold for $\ba_0=\bw^{(0)}$, then
we can conclude from the random initialization of $\bw^{(0)}$
that (\ref{le2eq1}), (\ref{le2eq2}), (\ref{le3eq1})
and (\ref{le3eq2}) hold with $\gamma_n^* = c_{21} \cdot (\log n)^2$ and
$B_n=c_{22} \cdot (\log n)^2$.
From this and Lemma \ref{le2}
and Lemma \ref{le3}
and the assumptions on $L_n$ and $t_n$ in Theorem \ref{th1}
we get that (\ref{le1eq1}) and (\ref{le1eq2})
hold.

In the {\it fifth step of the proof} we show
\[
\PROB(A_n^c) \leq \frac{c_{28}}{\beta_n^3}.
\]
To do this, we bound the probability that the weight
vector $\bw^{(0)}$ does not satisfy the first condition
in the definition of the event $A_n$ by considering a sequential
choice of the weights in the $K_n$ fully connected
neural networks which we compute in parallel.
In a single fully connected neural network the probability
that all $(L-1) \cdot r \cdot (r+1) + r \cdot (d+1)$ weights satisfy
the conditions of Lemma \ref{le7} for $j_1$ is bounded
from below by
\[
\left(
\frac{2 \cdot \log n}{40 d \cdot (\log n)^2}
\right)^{(L-1) \cdot r \cdot (r+1)}
\cdot
\left(
\frac{2 \cdot \log n}{2 n^{\tau}}
\right)^{r \cdot (d+1)}.
\]
Hence in the first $K_n/(K^2+1)^{3d}$ many fully connected networks
the probability that the condition for $j_1$ is never satisfied
is bounded from above by
\[
\left(
1 - \left(
\frac{1}{20 d \cdot \log n}
\right)^{(L-1) \cdot r \cdot (r+1)}
\cdot
\left(
\frac{\log n}{n^{\tau}}
\right)^{r \cdot (d+1)}
\right)^{K_n/(K^2+1)^{3d}}.
  \]
  This implies that all conditions of Lemma \ref{le7} are satisfied outside of
an event of probability
\[
(K^2+1)^{3d} \cdot
\left(
1 - \left(
\frac{1}{20 d \cdot \log n}
\right)^{(L-1) \cdot r \cdot (r+1)}
\cdot
\left(
\frac{\log n}{n^{\tau}}
\right)^{r \cdot (d+1)}
\right)^{K_n/(K^2+1)^{3d}},
\]
which is for large $n$ less than
\[
(K^2+1)^{3d} \cdot
\left(
1 -
\frac{1}{n^{r}}
\right)^{n^{r} \cdot \log n}
\leq
(K^2+1)^{3d} \cdot e^{- \log n} = \frac{(K^2+1)^{3d}}{n}
\leq
\frac{1}{2} \cdot \frac{c_{28}}{\beta_n^3}
\]
because of (\ref{th1eq2}) and $0 < \tau < 1/(d+1)$.
Furthermore, we can conclude from Markov's inequality
\[
\PROB \left\{
\frac{1}{n} \sum_{i=1}^n Y_i^2 > \beta_n^3
\right\}
\leq
\frac{
  \EXP \left\{
\frac{1}{n} \sum_{i=1}^n Y_i^2
  \right\}
}{
\beta_n^3
}
=
\frac{\EXP\{ Y^2\}}{\beta_n^3
}.
\]

In the {\it sixth step of the proof} we show
\[
\limsup_{n \rightarrow \infty} \EXP T_{3,n} \leq 0. 
\]
We have
\begin{eqnarray*}
  &&
  \frac{1}{
  1+\epsilon}
\cdot
 \EXP \left\{
T_{3,n} 
\right\}
\\
&&
\leq
\int_0^{4 \cdot \beta_n^2}
\PROB \Bigg\{
(\EXP\{ |m_n(X)-T_{\beta_n} Y|^2 \cdot 1_{[-\alpha_n,\alpha_n]^d}(X) | \D_n\}
\\
&&
\hspace*{1cm}
    -
    \frac{1}{n} \sum_{i=1}^n |m_n(X_i)- T_{\beta_n} Y_i|^2 \cdot 1_{[-\alpha_n,\alpha_n]^d}(X_i)) \cdot 1_{A_n}
    >t
    \Bigg\} \, dt
    \\
    &&
    \leq
    \frac{1}{n^{1/4}}
    +
    \int_{1/n^{1/4}}^{4 \cdot \beta_n^2}
\PROB \Bigg\{
    (\EXP\{ |m_n(X)-T_{\beta_n} Y|^2 \cdot 1_{[-\alpha_n,\alpha_n]^d}(X)| \D_n\}
\\
&&
\hspace*{2cm}
    -
    \frac{1}{n} \sum_{i=1}^n |m_n(X_i)- T_{\beta_n} Y_i|^2 \cdot 1_{[-\alpha_n,\alpha_n]^d}(X_i)) \cdot 1_{A_n}
    >t
    \Bigg\} \, dt.
\end{eqnarray*}

We want to derive a bound for the above probability. For this we can assume without loss of generality that $A_n$ holds. Due to the fourth step of
the proof it is then possible to apply Lemma \ref{le1} and to conclude
\[
\|\bw^{(t_n)}-\bw^{(0)}\|_\infty
  \leq
  \|\bw^{(t_n)}-\bw^{(0)}\|
    \leq c_{29} \cdot (\log n)^2. 
\]
From this and the choice of $\bw^{(0)}$ we can conclude that
$m_n$ is contained in the function space
\[
\{ T_{\beta_n} f \cdot 1_{[-\alpha_n,\alpha_n]^d} \, : \, f \in \F \}
\]
where 
$\F$ is defined as in Lemma \ref{le4}
with $C = c_{30} \cdot \sqrt{K_n} \cdot (\log n)^2$, $B=c_{31} \cdot (\log n)^2$
and $A = c_{32} \cdot n^{\tau}$. Application of Lemma \ref{le4}
together with standard bounds of empirical process theory (cf.,
Theorem 9.1 in Gy\"orfi et al. (2002)) yields
\begin{eqnarray*}
  &&
\PROB \Bigg\{
(\EXP\{ |m_n(X)-T_{\beta_n} Y|^2 \cdot 1_{[-\alpha_n,\alpha_n]^d}(X)| \D_n\}
\\
&&
\hspace*{3cm}
    -
    \frac{1}{n} \sum_{i=1}^n |m_n(X_i)- T_{\beta_n} Y_i|^2
    \cdot 1_{[-\alpha_n,\alpha_n]^d}(X_i)) \cdot 1_{A_n}
    >t
    \Bigg\}
    \\
    &&
    \leq
    8 \cdot \left(
c_{33}\cdot \frac{\beta_n}{t/8}
\right)^{
  c_{34} \cdot (\log n)^{c_{35}} \cdot n^{\tau \cdot d} \cdot \left(
\frac{c_{36} \cdot \sqrt{K_n} \cdot (\log n)^2}{t/8}
  \right)^{d/k} + c_{37}
}
\cdot
\exp \left(
- \frac{n \cdot t^2}{128 \beta_n^4}
\right).\\
  \end{eqnarray*}
By choosing $k$ large enough the right-hand side above is for $t>1/n^{1/4}$
bounded from above by
\[
c_{38}
\cdot
\exp \left(
- \frac{n \cdot t^2}{256 \cdot \beta_n^4}
\right)
\leq
c_{39}
\cdot
\exp \left(
- \frac{\sqrt{n}}{256 \cdot \beta_n^4}
\right).
\]
Consequently, we get
\[
 \EXP \left\{
T_{3,n} 
\right\}
\leq
(1+\epsilon) \cdot \left(
\frac{1}{n^{1/4}}+ 4 \beta_n^2 \cdot
c_{40}
\cdot
\exp \left(
- \frac{\sqrt{n}}{256 \beta_n^4}
\right)
\right)
\rightarrow 0
\quad (n \rightarrow \infty).
\]

In the {\it seventh step of the proof} we show
\[
\limsup_{n \rightarrow \infty} \EXP\{ T_{7,n} \} \leq 2 \epsilon.
\]
W.l.o.g. we assume $\|\bar{m}\|_\infty \leq \beta_n$. The boundedness
of $m_n$ by $\beta_n$ together with the fifth step imply
\begin{eqnarray*}
	\EXP\{ T_{7,n} \} &\leq&
	2 \cdot 4 \beta_n^2 \cdot \PROB(A_n^c) + 2 \int |\bar{m}(x)-m(x)|^2 \PROB_X (dx)
	\\
	&\leq&
	2 \cdot 4 \beta_n^2 \cdot \frac{c_{28}}{\beta_n^3}
	+ 2 \epsilon.
\end{eqnarray*}
Because of $\beta_n \rightarrow \infty$ $(n \rightarrow \infty)$
this implies the assertion of the seventh step.

In the {\it eighth step of the proof} we bound
\[
\EXP T_{6,n}.
\]
If $A_n$ holds, then we can apply Lemma \ref{le1}, which
together with Lemma \ref{le8} (which we can apply if we use
that the norm of the gradient of our nonlinear function
is larger than the sum of the squares of the partial derivatives
with respect to the weights corresponding only to the last layer $L$)
yields
\begin{eqnarray*}
  &&
  \frac{1}{n} \sum_{i=1}^n |f_{\bw^{(t_n)}}(X_i) - Y_i|^2 \cdot 1_{[-\alpha_n,\alpha_n]^d}(X_i)
  \\
  &&
  \leq
  \frac{1}{n} \sum_{i=1}^n |f_{\bw^{(t_n)}}(X_i) - Y_i|^2 \cdot 1_{[-\alpha_n,\alpha_n]^d}(X_i)
  + c_2 \cdot
  \sum_{k=1}^{K_n} ((\bw^{(t_n)})_{1,1,k}^{(L)})^2
  \\
  &&
  =
F_n(\bw^{(t_n)})
  \\
  &&
  \leq
  F_n(\bw^{(t_n-1)}) - \frac{1}{2 L_n} \cdot \| \nabla_\bw F_n(\bw^{(t_n-1)}) \|^2
  \\
  &&
  \leq
  F_n(\bw^{(t_n-1)}) - \frac{1}{2 L_n} \cdot 4 \cdot c_2 \cdot
  ( F_n(\bw^{(t_n-1)}) -  F_n((\bw^*)^{(t_n-1)}))
  \\
  &&
  =
  \left(
1 - \frac{2 \cdot c_2}{ L_n}  
\right)
\cdot
F_n(\bw^{(t_n-1)})
+
\frac{2 \cdot c_2}{ L_n} \cdot
F_n((\bw^*)^{(t_n-1)})
\\
&&
\leq
  \left(
1 - \frac{2 \cdot c_2}{ L_n}  
\right)^2
\cdot
F_n(\bw^{(t_n-2)})
+
\frac{2 \cdot c_2}{ L_n} \cdot
F_n((\bw^*)^{(t_n-1)})
\\
&&
\hspace*{3cm}
+
\frac{2 \cdot c_2}{ L_n} \cdot
  \left(
1 - \frac{2 \cdot c_2}{ L_n}  
\right)
F_n((\bw^*)^{(t_n-2)})
\\
&&
\leq \dots
\\
&&
\leq
  \left(
1 - \frac{2 \cdot c_2}{ L_n}  
\right)^{t_n}
\cdot
F_n(\bw^{(0)})
+
\sum_{k=1}^{t_n}
\frac{2 \cdot c_2}{ L_n} \cdot
  \left(
1 - \frac{2 \cdot c_2}{ L_n}  
\right)^{k-1}
F_n((\bw^*)^{(t_n-k)}).
\end{eqnarray*}
This implies
\begin{eqnarray*}
  &&
  \EXP \left\{
T_{6,n} 
\right\}
\\
&&
\leq
(1+\epsilon)^2
\cdot
  \left(
1 - \frac{2 \cdot c_2}{ L_n}  
\right)^{t_n}
\cdot \EXP\{Y^2\}
\\
&&
\quad +
(1+\epsilon)^2
\cdot
\sum_{k=1}^{t_n}
\frac{2 \cdot c_2}{ L_n} \cdot
  \left(
1 - \frac{2 \cdot c_2}{ L_n}  
\right)^{k-1}
\Big(
\EXP\{ F_n((\bw^*)^{(t_n-k)}) \cdot 1_{A_n}\}
\\
&&
\hspace*{8cm}
-
   \EXP\{ |m(X)-Y|^2 \} \cdot \PROB(A_n)
   \Big)
   \\
   &&
   \quad
   + ((1+\epsilon)^2 - 1) \cdot
   \EXP\{ |m(X)-Y|^2 \}.
  \end{eqnarray*}
We have
\[
  \left(
1 - \frac{2 \cdot c_2}{ L_n}  
\right)^{t_n}
\leq
\exp \left(
\frac{-2 \cdot c_2 \cdot t_n}{L_n}
\right)
\leq
\exp
\left(
- c_{41} \cdot \log n
\right)
\rightarrow 0
\quad (n \rightarrow \infty)
\]
and
\begin{eqnarray*}
	&&
	\EXP\{ F_n((\bw^*)^{(t_n-k)}) \cdot 1_{A_n} \}
	-
	\EXP\{ |m(X)-Y|^2 \} \cdot \PROB(A_n)
        \\
        &&
        =
        \EXP\Big\{\frac{1}{n} \sum_{i=1}^{n}|f_{((\bw^*)^{(t_n-k)})}(X_i)-Y_i|^2 \cdot 1_{[-\alpha_n,\alpha_n]^d}(X_i) \cdot 1_{A_n}
        +
        c_2 \cdot \sum_{j=1}^{K_n}|(\bw^*)^{(L)}_{1,1,j}|^2
        \Big\}
        \\
        &&
        \quad
        	-
	        \EXP\{ |m(X)-Y|^2 \} \cdot \PROB(A_n)
                \\
                &&
                = \Bigg( \EXP\Big\{\frac{1}{n} \sum_{i=1}^{n}|f_{((\bw^*)^{(t_n-k)})}(X_i)-Y_i|^2 \cdot 1_{[-\alpha_n,\alpha_n]^d}(X_i) \cdot 1_{A_n}
                \Big\}
                \\
                &&
                \hspace*{1cm}
                - (1+\epsilon) \cdot
                \EXP\Big\{\frac{1}{n} \sum_{i=1}^{n}|f_{((\bw^*)^{(t_n-k)})}(X_i)-
                T_{\beta_n} Y_i|^2 \cdot 1_{[-\alpha_n,\alpha_n]^d}(X_i) \cdot 1_{A_n}
                \Big\} \Bigg)
                \\
                &&
                \quad +
                (1+\epsilon) \cdot
                \Bigg(
                \EXP\Big\{\frac{1}{n} \sum_{i=1}^{n}|f_{((\bw^*)^{(t_n-k)})}(X_i)-
                T_{\beta_n} Y_i|^2 \cdot 1_{[-\alpha_n,\alpha_n]^d}(X_i) \cdot 1_{A_n}
                \Big\}
                \\
                &&
                \hspace*{1cm}
                -
                \EXP \Big\{ \EXP\Big\{|f_{((\bw^*)^{(t_n-k)})}(X)-
                T_{\beta_n} Y|^2 \cdot 1_{[-\alpha_n,\alpha_n]^d}(X) \cdot 1_{A_n} \Big| \D_n
                \Big\} \Big\}
                \Bigg)
                \\
                &&
                \quad
                +
                (1+\epsilon) \cdot \Bigg(
                \EXP \Big\{ |f_{((\bw^*)^{(t_n-k)})}(X)-
                T_{\beta_n} Y|^2 \cdot 1_{[-\alpha_n,\alpha_n]^d}(X) \cdot 1_{A_n} \Big\}
                \\
                &&
                \hspace*{1cm}
                - (1+\epsilon)^2 \cdot
         \EXP \Big\{ |f_{((\bw^*)^{(t_n-k)})}(X)-
                Y|^2 \cdot 1_{[-\alpha_n,\alpha_n]^d}(X) \cdot 1_{A_n} \Big\}       \Bigg)
                \\
                &&
                \quad
                +
                (1+\epsilon)^2 \cdot \Bigg(
         \EXP \Big\{ |f_{((\bw^*)^{(t_n-k)})}(X)-
         Y|^2 \cdot
         1_{[-\alpha_n,\alpha_n]^d}(X) \cdot 1_{A_n} \Big\}
         \\
         &&
         \hspace*{8cm}
         -
                \EXP\{ |m(X)-Y|^2 \} \cdot \PROB(A_n)
                \Bigg)
                \\
                &&
                \quad
                + ((1+\epsilon)^2 -1) \cdot  \EXP\{ |m(X)-Y|^2 \}
                + c_2 \cdot \sum_{j=1}^{K_n}|(\bw^*)^{(L)}_{1,1,j}|^2
                	\\
	                && = T_{8,n}+T_{9,n}+T_{10,n}+T_{11,n}+T_{12,n}+T_{13,n}.
                        \end{eqnarray*}
Similar to the second step, we obtain
	\begin{align*}
	  \limsup_{n \rightarrow \infty }T_{8,n} \leq 0 \quad \mbox{and}
          \quad
          \limsup_{n \rightarrow \infty } T_{10,n} \leq 0.
	\end{align*}
According to the assertion of Lemma \ref{le6} we know that $f_{\bw^*}$ is bounded.
Thus, we get as in the sixth step
\begin{align*}
	\limsup_{n \rightarrow \infty} T_{9,n} \leq 0.
\end{align*}
The choice of $\bar{m}$ and Lemma \ref{le7} imply
\begin{align*}
  &  T_{11,n}/(1+\epsilon)^2 \\
  &\leq \EXP \Big\{ \int |f_{((\bw^*)^{(t_n-k)})}(x)-m(x)|^2 \PROB_X(dx) \cdot 1_{A_n}
  \Big\} \\
  &\leq 2 \EXP \Big\{ \int |f_{((\bw^*)^{(t_n-k)})}(x)-\bar{m}(x)|^2
  \PROB_X(dx) \cdot 1_{A_n} \Big\} + 2\int |\bar{m}(x)-m(x)|^2 \PROB_X(dx)\\
  	& \leq  c_{42} \cdot \left( \frac{1}{K} + \frac{K^{12d}}{n^2} 
	+\left(\frac{K^{6d}}{n} +1\right)^2 \PROB_X(\mathbb{R}^d\setminus [-K,K]^d)\right)+ 2 \epsilon,
\end{align*}
from which we can conclude
\begin{align*}
	\limsup_{n \rightarrow \infty} T_{11,n} \leq c_{42} \cdot \left( \frac{1}{K} +  \PROB_X(\mathbb{R}^d\setminus [-K,K]^d)\right)+ 2\epsilon \cdot (1+\epsilon)^2.
\end{align*}
Furthermore we obtain by Lemma \ref{le7}
\begin{align*}
	T_{13,n} \leq c_{43}\cdot (K^2+1)^{3d} \cdot \left( \frac{1}{(K^2+1)^{2d}} \right)^2.
\end{align*}
Summarizing the above results yields
\begin{eqnarray*}
  &&
\limsup_{n \rightarrow \infty} \sum_{k=1}^{t_n}
\frac{2 \cdot c_2}{ L_n} \cdot
  \left(
1 - \frac{2 \cdot c_2}{ L_n}  
\right)^{(k-1)}
\Big(
\EXP\{ F_n((\bw^*)^{(t_n-k)}) \cdot 1_{A_n} \}
\\
&&
\hspace*{8cm}
-
   \EXP\{ |m(X)-Y|^2 \} \cdot \PROB(A_n)
   \Big)
   \\
   &&
  \leq
  c_{44} \cdot \left( \frac{1}{K} +  \PROB_X(\mathbb{R}^d\setminus [-K,K]^d)
+ \epsilon
  \right)+ c_{43}\cdot (K^2+1)^{3d} \cdot \left( \frac{1}{(K^2+1)^{2d}} \right)^2
\end{eqnarray*}
and 
\begin{align*}
\limsup_{n \rightarrow \infty}
  \EXP \left\{
T_{6,n} 
\right\}
&\leq
(1+\epsilon)^2 \left( c_{44} \cdot \left(
 \frac{1}{K}  + \PROB_X \left(
 \Rd \setminus [-K,K]^d
 \right)
+ \epsilon
 \right)
  + c_{45} \cdot \frac{1}{(K^2+1)^d} \right)
  \\
  &\qquad +
((1+\epsilon)^2 - 1) \cdot
   \EXP\{ |m(X)-Y|^2 \}.
\end{align*}

In the {\it ninth step of the proof} we finish the proof
of Theorem \ref{th1}. The results of steps 1,2,3,6,7 and 8 imply
for $K \rightarrow \infty$
\[
\limsup_{n \rightarrow \infty}
\EXP     \int | m_n(x)-m(x)|^2 \PROB_X (dx)
\leq
c_{46} \cdot \epsilon.
\]
With $\epsilon \rightarrow 0$ we get the assertion.
\hfill $\Box$

\end{document}